\documentclass{article}%
\usepackage{amsthm,amssymb}
\catcode`\@=11
\@addtoreset{equation}{section}

\newtheorem{thm}{Theorem}[section]
\newtheorem{prop}{Proposition}[section]

\newtheorem{rem}{Remark}[section]
\newtheorem{lem}{Lemma}[section]
\newtheorem{dfn}{Definition}[section]

\newcommand{\eps}{\varepsilon}

\newcommand{\codim}{{\mathrm{codim}\,}}

\newcommand{\Lip}{{\mathrm{Lip}}}
\newcommand{\Crit}{{\mathrm{Crit}}}

\newcommand{\R}{{\mathbb{R}}}

\newcommand{\kbf}{{\mathbf{k}}}

\newcommand{\xbf}{{\mathbf{x}}}

\newcommand{\zbf}{{\mathbf{z}}}

\newcommand{\sbf}{{\mathbf{s}}}
\newcommand{\qbf}{{\mathbf{q}}}
\newcommand{\ubf}{{\mathbf{u}}}
\newcommand{\vbf}{{\mathbf{v}}}

\newcommand{\A}{{\mathbb{A}}}
\newcommand{\C}{{\mathbb{C}}}
\newcommand{\Z}{{\mathbb{Z}}}
\newcommand{\T}{{\mathbb{T}}}

\newcommand{\calA}{\mathcal{A}}
\newcommand{\calB}{\mathcal{B}}
\newcommand{\calC}{\mathcal{C}}
\newcommand{\calD}{\mathcal{D}}

\newcommand{\calH}{\mathcal{H}}
\newcommand{\calK}{\mathcal{K}}
\newcommand{\calL}{\mathcal{L}}
\newcommand{\calM}{\mathcal{M}}

\newcommand{\calF}{\mathcal{F}}
\newcommand{\calP}{\mathcal{P}}
\newcommand{\calQ}{\mathcal{Q}}
\newcommand{\calR}{\mathcal{R}}

\newcommand{\calT}{\mathcal{T}}

\newcommand{\calX}{\mathcal{X}}

\title{Degenerate billiards in celestial mechanics
\author{Sergey Bolotin\thanks{Supported by the RFBR grant of the Russian Academy of Sciences "Modern problems of classical dynamics"   
(project 15-01-03747a)}
\\
Department of Mathematics\\
University of Wisconsin--Madison\\
and\\
Moscow Steklov Mathematical Institute
} }

\begin{document}

\maketitle

\begin{abstract}
In an ordinary billiard  trajectories of a Hamiltonian system are
elastically reflected after a collision with a hypersurface (scatterer).
If the scatterer is a submanifold of codimension more than one,
we say that the billiard is degenerate.
Degenerate billiards appear as limits of systems  with singularities in celestial mechanics.
We prove the existence of trajectories  of such systems
shadowing trajectories of the corresponding degenerate billiards.
This research is motivated by the problem of second species solutions of Poincar\'e.
\end{abstract}

\section{Introduction}
\subsection{Degenerate billiards}

\label{sec:bill}

Consider  a Hamiltonian system $(M,H)$ with the configuration space $M$ and a classical
smooth\footnote{$C^4$ is enough. We do not attempt to lower regularity since in applications
to celestial mechanics $H$ is real analytic.}
Hamiltonian $H$ on the  phase space  $T^*M$:
\begin{equation}
\label{eq:H}
H(q,p)=\frac12\|p-w(q)\|^2+ W(q),
\end{equation}
Here  $\|\,\cdot\,\|$ is a Riemannian metric on $M$, and $w$ a covector field representing gyroscopic (or magnetic) forces.
The symplectic structure $dp\wedge dq$ on $T^*M$  is standard, so we do not include it in the notation.
Let\footnote{We use the same notation $\|\,\cdot\, \|$ for the norm of a vector and a covector.}
\begin{equation}
\label{eq:L}
L(q,\dot q)=\max_p (\langle p,\dot q\rangle-H(q,p))=\frac12\|\dot q\|^2+\langle w(q),\dot q\rangle -W(q)
\end{equation}
be the corresponding Lagrangian.

\begin{rem}
A transformation $p\to p+w(q)$  replaces $H$ with a natural Hamiltonian
\begin{equation}
\label{eq:natural}
H(q,p)=\frac12\|p\|^2+ W(q).
\end{equation}
Then the symplectic structure  is twisted
\begin{equation}
\label{eq:omega}
\omega=dp\wedge dq+\pi^*\Omega,\qquad \pi:T^*M\to M,
\end{equation}
where  $\Omega=dw$ is the gyroscopic 2-form\footnote{The covector field $w$ is regarded as a 1-form on $M$.}
on $M$.

Conversely, if the 2-form $\Omega$  is exact, we can  make the symplectic structure (\ref{eq:omega}) standard and the Hamiltonian
takes the form  (\ref{eq:H}).
The twisted symplectic structure is   convenient for many purposes,
i.e.\ reduction of symmetry. However for simplicity
we will use the  standard form $dp\wedge dq$ and the Hamiltonian (\ref{eq:H}).
\end{rem}

Let $N\subset M$ be a submanifold in $M$ which is called a scatterer.  Suppose that when a
trajectory\footnote{A trajectory $q(t)$ is the projection of a solution $(q(t),p(t))$ to the configuration space.}
$q(t)$  meets the scatterer at a collision point
$x=q(\tau)\in N$, it is reflected according to the elastic reflection law\footnote{Here $p_-$ is the momentum after the collision,
and $p_+$ before the collision. The strange notation is chosen to fit with the notation for the initial and final momenta of a collision orbit $\gamma:[t_-,t_+]\to M$ in (\ref{eq:coll}).}
\begin{eqnarray}
\label{eq:Deltap} \Delta p(\tau)=p_--p_+\perp T_x N,\qquad p_\pm=p(\tau\mp 0),\\
\Delta H(\tau)=H(x,p_-)-H(x,p_+)=0.
\label{eq:DeltaH}
\end{eqnarray}
Thus the tangent component $y\in T_x^*N$ of the momentum $p\in T_x^*M$ and the energy $H=E$ are preserved.
We always assume that the momentum has a jump at the collision: $\Delta p(\tau)\ne 0$.
Then also collision velocities
$v_\pm=\dot q(\tau\mp 0)$
have a jump $\Delta v(\tau)=v_--v_+\ne 0$ orthogonal to $N$ with respect to the Riemannian metric and they are not tangent to the scatterer: $v_\pm\notin T_xN$. By conservation of energy
\begin{equation}
\label{eq:energy}
\|v_+\|^2=\|v_-\|^2=2(E-W(x)).
\end{equation}

When $N$ is a hypersurface bounding a domain $\Omega$ in $M$,   we obtain
a usual billiard system $(\Omega,N,H)$. If
$$
d=\codim N>1,
$$
we say that $(M,N,H)$ is a  degenerate billiard. We do not assume $N$ to be connected, it may have connected components
of different dimension.

Trajectories of the  degenerate billiard having collisions with
$N$ form a zero measure set  in the phase space.
Moreover  $p_+$ does not determine $p_-$ uniquely by
(\ref{eq:Deltap})--(\ref{eq:DeltaH}): for given $p_+$ the set of possible $p_-$
has dimension $d-1$.
Thus the past of a collision trajectory does not determine its future.
The simplest case is  when $N$ is a discrete set in $M$, then only condition (\ref{eq:energy}) remains
and a trajectory can be reflected in any direction.

We are interested in trajectories  $\gamma:[\alpha,\beta]\to M$ with multiple collisions
which are called collision chains.  They  are extremals of the
action functional
\begin{equation}
\label{eq:I}
I(\gamma)=\int_\alpha^\beta L(\gamma(t),\dot\gamma(t))\,dt=\sum_{j=0}^nI(\gamma_j),\qquad \gamma_j=\gamma|_{[t_j,t_{j+1}]},
\end{equation}
on the set of curves $\gamma:[\alpha,\beta]\to M$ with fixed end points $a=\gamma(\alpha)$ and $b=\gamma(\beta)$, subject to the
constraints $\gamma(t_j)= x_j\in N$, $j=1,\dots,n$, for some sequence
$\alpha=t_0<t_1<\dots <t_n<t_{n+1}=\beta$. Here $x_j$, $t_j$ are independent variables.
Each segment $\gamma_j=\gamma|_{[t_j,t_{j+1}]}$ is a collision orbit joining points in $N$ and
\begin{equation}
\label{eq:elastic2}
\Delta p(t_j)\perp T_{x_j} N,\quad\Delta H(t_j)=0.
\end{equation}
 We also require the  jump condition
 \begin{equation}
 \label{eq:jump}
 \Delta p(t_j)\ne 0.
 \end{equation}

An infinite collision chain $\gamma:\R\to M$
is a concatenation  of a sequence $\gamma=(\gamma_j)_{j\in\Z}$ of
collision orbits $\gamma_j:[t_j,t_{j+1}]\to M$    such that the
elastic reflection law  (\ref{eq:elastic2})--(\ref{eq:jump})  is satisfied at each
collision.

An evident source of degenerate billiards are billiards with thin scatterers.
 Let $N$ be a submanifold in $M$   and $N_\eps $
its tubular $\eps$-neighborhood with the boundary $\Sigma_\eps=\partial N_\eps$.
Consider the billiard system $(\Omega_\eps,\Sigma_\eps,H)$ in the domain $\Omega_\eps=M\setminus N_\eps$ with the boundary
 $\partial \Omega_\eps=\Sigma_\eps$ and Hamiltonian $H$. As $\eps\to 0$, it  approaches the degenerate billiard
$(M,N,H)$ with the scatterer $N$. In \cite{Bol:MIAN} it is proved that for small $\eps>0$ nondegenerate collision chains of this
degenerate billiard are shadowed by trajectories of the billiard system in $\Omega_\eps$.
For a discrete set $N$, this was shown earlier in \cite{Chen}, see also \cite{Bol-Tre:ANTI}.

The goal of the present paper is to show how degenerate billiards appear in Hamiltonian systems
with Newtonian singularities. The motivation is the study of periodic and chaotic second species solutions of Poincar\'e
in celestial mechanics \cite{Poincare}, see  section \ref{sec:examples}. It turns out that the problem is reduced to understanding the corresponding degenerate billiard.

The results of this paper generalize some results of \cite{Bol-Mac:centers} where $N$ was  a discrete set
and of \cite{Bol:DCDS,Bol-Neg:DCDS} where $N$ was 2-dimensional.

\subsection{Systems with Newtonian singularities}
\label{sec:sing}

Consider a Hamiltonian system $(M\setminus N,H_\mu)$ on $T^*(M\setminus N)$ with a classical smooth\footnote{$C^4$ is enough.}
Hamiltonian
\begin{equation}
\label{eq:Hmu}
H_\mu(q,p)=\frac12\|p-w_\mu(q)\|_\mu^2+W_\mu(q)+\mu V(q,\mu)
\end{equation}
depending on a small parameter $\mu\in (-\mu_0,\mu_0)$. Here $\|\,\cdot\,\|_\mu$ is a Riemannian metric on $M$,
smoothly depending on   $\mu$,
and $w_\mu$ and $W_\mu$ are covector field and a function on $M$ smoothly depending on $\mu$.
The potential $V$ is smooth on $M\setminus N$ but undefined on $N$.

We say that $V$ has a Newtonian singularity on $N$
if in a tubular neighborhood of $N$
there exists a smooth positive function $\phi$  such that
\begin{equation}
\label{eq:V}
V(q,\mu)=-\frac{\phi(q,\mu)}{d_\mu(q,N)}.
\end{equation}
The distance $d_\mu$ is defined by the Riemannian metric $\|\,\cdot\,\|_\mu$.  If $\mu<0$, the singular force  is repelling (like the Coulomb force),
and if $\mu>0$ attracting (like the gravitational force).

For $\mu=0$  the singularity disappears  and we obtain the Hamiltonian system $(M,H_0)$ with Hamiltonian $H_0=H$
as in  (\ref{eq:H}).
The perturbation consists of two parts: regular perturbation
which is a smooth function on $T^*M$, and a singular part $\mu V$.

We are interested in nearly collision trajectories  of system $(M\setminus N,H_\mu)$
which pass $O(\mu)$-close to $N$.
Their  limits  as $\mu\to 0$ are collision chains of the degenerate billiard $(M,N,H_0)$
with Hamiltonian $H_0$ and scatterer $N$.
We will give precise statements in section \ref{sec:main}.

\begin{rem}
For $\mu<0$ (repelling force) trajectories of system $(M\setminus N,H_\mu)$ do not have collisions with $N$.
For $\mu>0$ collisions may appear.
However, the Hamiltonian flow on the energy level $\{H_\mu=E\}$ is regularizable (see section \ref{sec:reg}),
since collisions with $N$  are of the type of double collisions in celestial mechanics.
After a change of variables and a time reparametrization we obtain  a smooth flow without singularities.
\end{rem}

\subsection{Examples}

\label{sec:examples}

1. The $n$ center problem. Suppose a particle moves  in $\R^3$ under the gravitational forces of $n$ fixed centers $a_1,\dots,a_n$
with small masses $m_i=\mu \alpha_i$, $0<\mu\ll 1$. By a time change  $t\to t/\sqrt{\mu}$ this is equivalent to the case of
centers of finite masses $\alpha_i$ and large energy  of order $\mu^{-1}$ of the particle. Then
\begin{equation}
\label{eq:ncent}
H_\mu(q,p)=\frac12|p|^2+\mu V(q),\qquad V(q)=-\sum_{i=1}^n\frac{\alpha_i}{|q-a_i|},\qquad q\in\R^3 .
\end{equation}
The limit system is the degenerate billiard $(\R^3,\{a_1,\dots,a_n\}, H_0)$ with a finite scatterer and Hamiltonian $H_0=|p|^2/2$.
Collision chains are polygons with vertices $a_i$. If $n\ge 4$ there is a Cantor set of collision chains, see \cite{Knauf}.
In fact  $n$ center problem in $\R^3$ has chaotic invariant sets on positive energy levels for $n\ge 3$ and any $\mu>0$ for purely topological reasons, 
see \cite{Bol-Neg:reg}.

2. A more  realistic example is the restricted $n+2$ body problem.
Then the bodies $a_1,\dots,a_n$ with small masses $m_i=\mu\alpha_i$
move around the  Sun with mass $1-\mu$ along circular orbits with the same angular velocity $\omega\in\R^3$.
An Asteroid of negligible mass moves under the action of the gravitational  forces of
the Sun and the small bodies.
Then in a rotating coordinate frame,
$$
H_\mu(q,p)=\frac12|p-\omega\times q|^2-\frac{1}{|q|}+\mu V(q)+O(\mu),
$$
where $V(q)$ is as in (\ref{eq:ncent}). The corresponding degenerate billiard $(\R^3\setminus \{0\}, \{a_1,\dots,a_n\}, H_0)$
has the same scatterer as in   example 1, but now $H_0$ is the Hamiltonian of the Kepler problem in
a rotating coordinate frame. Because of this the set of collision chains with fixed energy $H_0=E$ (called Jacobi integral) is very rich: 
already for $n=1$  it is a Cantor set, and there is a hyperbolic chaotic set of  shadowing orbits, see \cite{Bol-Mac:centers}.
Shadowing periodic  orbits are called second species solutions of Poincar\'e.
They are well studied for the circular restricted 3 body problem, see e.g.\ \cite{Gomez,Perko,Mar-Nid,Bol-Mac:centers,FNS}
and for the elliptic restricted 3 body problem, see e.g.\ \cite{Gomez,Bol:ELL,Bol:DCDS}.
Poincar\'e \cite{Poincare} considered  the nonrestricted 3 body problem, see example 4 below.

3. The $n$ body problem with small masses $m_i=\mu \alpha_i$  (or finite masses and large energy). Then
after a time change,
$$
H_\mu(q,p)= \frac12\|p\|^2+\mu V(q),\qquad q=(q_1,\dots,q_n)\in\R^{3n},
$$
where
$$
\|p\|^2=\sum_{i=1}^n\frac{|p_i|^2}{\alpha_i},  \quad V(q)=\sum_{i\ne j}\frac{\alpha_i\alpha_j}{|q_i-q_j|}.
$$
The  limit system is the degenerate billiard $(\R^{3n},\Delta,H_0\}$ with the scatterer
$$
\Delta=\cup_{i\ne j}\{q\in\R^{3n}:q_i=q_j\}
$$
and Hamiltonian $H_0=\|p\|^2/2$. The scatterer is not a manifold, so to obtain a billiard of the type studied in this paper 
we need to exclude from $\Delta$  multiple collisions. Dynamics of this billiard is  finite: after a bounded number of collisions
the bodies escape to infinity. This is a deep result proved  in \cite{Burago}, see also \cite{Mont}.
Hence this degenerate billiard does not have invariant sets, in particular  it has no periodic orbits and no chaotic hyperbolic
sets which are of interest to us.

4. The most important example was introduced by  Poincar\'e \cite{Poincare}.
Consider the  $n+1$ body problem  with one of the masses $m_{0}$
much larger than the rest. Set
$$
\frac{m_i}{m_{0}}=\mu\alpha_i,\qquad \sum_{i=1}^{n}\alpha_i=1,\qquad \mu\ll 1.
$$
We may assume that the center of mass is at rest: $\sum_{i=0}^{n}p_i=0$. Let  $q_i$
be the relative position of $m_i$ with respect to $m_{0}$.
Then after a time change we obtain the Hamiltonian
\begin{equation}
\label{eq:n_body}
H_\mu(q,p)=H_0(q,p) +\frac{\mu }2\bigg|\sum_{i=1}^{n}p_i\bigg|^2 -
\mu\sum_{i\ne j}\frac{\alpha_i\alpha_j}{|q_i-q_j|},
\end{equation}
where $q=(q_1,\dots,q_{n})\in\R^{3n}$ and
\begin{equation}
\label{eq:H0}
H_0=\sum_{i=1}^{n}\left(\frac{|p_i|^2}{2 \alpha_i}-\frac{\alpha_i}{|q_i| }\right).
\end{equation}
The Hamiltonian $H_0$ describes $n$ uncoupled Kepler problems which are of course integrable.
However for $n\ge 2$ the corresponding degenerate billiard $((\R^{3}\setminus\{0\})^{n},\Delta,H_0)$ has complicated chaotic dynamics.
Orbits of the $n+1$ body problem shadowing collision chains of  this billiard
are called second species solutions of Poincar\'e \cite{Poincare}.  Poincar\'e  discussed such solutions for the 3 body problem,
but did not provide a rigorous proof of their existence. There are many works of Astronomers on the subject but few mathematical results
(except for the restricted circular  3 body problem, see e.g.\ \cite{Perko,Gomez,Mar-Nid,Bol-Mac:centers} and
the elliptic restricted 3 body problem, see e.g.\ \cite{Gomez,Bol:ELL, Bol:DCDS}).
Some rigorous results for the unrestricted plane 3 body problem were proved in \cite{Bol-Neg:DCDS,Bol-Neg:RCD}
and for the 2 center - 2 body problem in \cite{Dolgop}.

\section{ Shadowing collision chains}

\label{sec:main}

\subsection{Discrete Lagrangian system of a degenerate billiard}

Before formulating   the  main results
we need to recall some definitions  from \cite{Bol:MIAN}, see also \cite{Bol-Neg:DCDS}.

The Hamiltonian is constant along collision chains of a degenerate billiard, so let us fix energy $H=E$.
The restriction of the Hamiltonian system $(M,H)$ to the energy level  will be denoted $(M,H=E)$.
 Trajectories $\gamma:[\alpha,\beta]\to M$  with energy $E$  are extremals of the Maupertuis action $J=J_E$:
\begin{equation}
\label{eq:J}
J(\gamma)=\int_\alpha^\beta g_E(\gamma(t),\dot\gamma(t))\,dt,
\end{equation}
i.e.\ geodesics\footnote{We identify curves which differ by an orientation preserving reparametrization.} of the Jacobi metric \cite{Arnold,AKN}
\begin{equation}
\label{eq:metric}
 g_E(q,\dot  q)=\max_p \{\langle p,\dot q\rangle:H(q,p)=E\}
 =\sqrt{2(E-W(q))}\|\dot q\|+\langle w(q),\dot q\rangle
\end{equation}
in the domain of possible motion
\begin{equation}
\label{eq:DE}
\calD_E =\{q\in M:W(q)<E\}.
\end{equation}

\begin{rem}
The metric $g_E$ is positive definite in the domain
$$
\{q\in M:W(q)+\|w(q)\|^2/2<E\},
$$
but not in $\calD_E$, so $g_E$ is not a Finsler metric in $\calD_E$.
However $g_E$ is convex in the velocity, so local calculus of variations works. In particular,
for any $x_0\in\calD_E$ there is $r>0$ such that a pair of points in the ball $B_r(x_0)$
is joined  by a geodesic in $B_r(x_0)$.
\end{rem}

For trajectories $\gamma$ with energy $E$,
$$
J(\gamma)=\int_\gamma p\, dq.
$$

When the energy is fixed, we denote the degenerate  billiard by $(M,N,H=E)$.
As in section 1, we call a  trajectory  $\gamma:[t_-,t_+]\to M$   a collision orbit if
its end points   lie in $N$ and there is no tangency and no early collisions with the scatterer:
\begin{equation}
\label{eq:coll}
\gamma(t_\pm)=a_\pm\in N,\qquad
v_\pm=\dot\gamma(t_\pm)\notin T_{a_\pm}N,\qquad
\gamma(t)\notin N,\quad t_-<t<t_+.
\end{equation}
In particular, $a_\pm\in\calD_E\cap N$.

We call $\gamma$ nondegenerate if it is nondegenerate as a critical point of $J$, i.e.\ the points $a_-$ and $a_+$ are non-conjugate.
Then there exist neighborhoods  $U_\pm\subset M$ of $a_\pm$ such that for all
$q_\pm\in U_\pm$ there exists an   orbit $\gamma(q_-,q_+)$ with energy $E$ joining $q_-$ and $q_+$, and it smoothly depends on $q_-,q_+$.
 The Maupertuis action
\begin{equation}
\label{eq:S}
S(q_-,q_+)=J(\gamma(q_-,q_+))
\end{equation}
is a smooth function on $U_-\times U_+$. The initial and final momenta of the orbit $\gamma$  are
$$
p_-=-D_{q_-}S,\quad p_+=D_{ q_+}S.
$$
The twist of the action function  is the linear transformation
$$
B(q_-,q_+)= D_{ q_-}D_{ q_+}S:T_{q_-}M\to T_{q_+}^*M,
$$
i.e.\ a bilinear form    on $T_{q_-}M\times T_{q_+}M$. Since the Hamiltonian system is autonomous, it is always degenerate:
\begin{equation}
\label{eq:twist}
B(a_-,a_+)v_-=0,\quad B^*(a_-,a_+)v_+=0 .
\end{equation}

 We say that the collision orbit $\gamma$ has nondegenerate
twist if the restriction
of the bilinear form $B(a_-,a_+)$ to $T_{a_-}N\times T_{a_+}N$    is nondegenerate.
For  this it is necessary that $v_\pm\notin
T_{a_\pm}N$, i.e.\ the collision orbit is not tangent to the scatterer $N$ at the end points.
For an ordinary billiard, when $N$ is a hypersurface, this is also sufficient for the nondegenerate twist,
but in general not for a degenerate billiard.

If $\gamma$ has nondegenerate  twist, the restriction of $S$   to a neighborhood
of $(a_-,a_+)$ in $N\times N$ is the generating function of a locally
defined symplectic map $f:V^-\to V^+$ of open sets
$V^\pm\subset T^*N$:
$$
f(x_{-},y_{-})=(x_{+},y_{+})\quad \Leftrightarrow\quad
y_+= D_{ x_+}S,\quad y_-=-D_{ x_-}S.
$$
Here $y_{\pm}=p_\pm|_{T_xN}\in T_{x_\pm}^*N$ are the tangent  projections  of the
collision momenta.
Hence $V^\pm\subset \calM_E$, where
\begin{eqnarray}
\label{eq:calM}
\calM_E=\{(x,y)\in T^*N: F(x,y)<E\},\label{eq:calM}\\
F(x,y)=\min_{p|_{T_xN}=y} H(x,p)=\frac12 \|y-a(x)\|^2+W(x).
\label{eq:F}
\end{eqnarray}
Here $a(x)=w(x)|_{T_xN}\in T_x^*N$.
The Riemannian metric is the induced metric on $N$.
Thus $F$ is the Hamiltonian on $T^*N$ corresponding to the Lagrangian $L|_{TN}$.

\begin{rem}
If the symplectic structure  (\ref{eq:omega}) is twisted, then,  locally, $\Omega=dw$,
where $w$ is   defined up to adding a differential $d\varphi$.
Then the  generating function   $S(x_-,x_+)$ is  defined up to adding a cocycle $\varphi(x_+)-\varphi(x_-)$.
\end{rem}

In general there may exist several (or none) nondegenerate
collision orbits  with   energy $E$ joining a
pair of points in $N$. Thus we obtain a collection   $\calL=\{ L_k\}_{k\in K}$ of action functions (\ref{eq:S}) on
open sets $U_k\subset N\times N$. Under the twist condition,    $L_k$ generates a
local symplectic map $f_k:V_k^-\to V_k^+$ of open sets in $\calM_E$.
We call the partly defined multivalued ``map''  $\calF=\{f_k\}_{k\in K}$ of $\calM_E$ the collision map,
or the scattering map of the degenerate billiard. It is analogous to the scattering map
of a normally hyperbolic invariant manifold, see \cite{Delshams}.
The degenerate billiard defines a discrete dynamical system    -- the skew product of the maps
$\calF=\{f_k\}_{k\in K}$ which is a map of a subset in $ K^\Z\times\calM_E$.

\begin{rem}
Computation of the collision map is usually difficult. See e.g.\ \cite{Bol:ELL,Bol-Neg:DCDS}
for the degenerate billiards appearing in the elliptic restricted 3 body problem and in the nonrestricted
plane 3 body problem.
\end{rem}

An orbit of $\calF$ is a pair $(\kbf,\zbf)$ of sequences $\kbf=(k_j)$,
$\zbf=(z_j)$, where $z_j=(x_j,y_j) \in
V_{k_j}^-\cap V_{k_{j-1}}^+$, such that $z_{j+1}=f_{k_j}(z_j)$.  The orbit
$(\kbf,\zbf)$ defines  a chain  of collision orbits
$\gamma_j$ joining $x_j$ with $x_{j+1}$.  The tangent collision  momenta of the collision chain are
\begin{equation}
\label{eq:y_j}
y_j=  D_{x_j} L_{k_{j-1}}(x_{j-1},x_{j})=-D_{x_j} L_{k_j}(x_j,x_{j+1}).
\end{equation}

Also without the twist condition,
the degenerate billiard $(M,N,H=E)$ can be viewed as a discrete Lagrangian
system (DLS) with multivalued Lagrangian $\calL=\{ L_k\}_{k\in K}$, see \cite{Bol-Tre:ANTI}.
Infinite collision chains  correspond to  critical points
$\xbf=(x_j)_{j\in\Z}$ of the  discrete action functional
\begin{equation}
\label{eq:Ak}
\calA_\kbf(\xbf)=\sum_{j\in\Z}  L_{k_j}(x_j,x_{j+1}),\qquad (x_j,x_{j+1})\in U_{k_j}.
\end{equation}
 For infinite collision chains, the   sum  makes no sense, so $A_\kbf(\xbf)$ is a
formal functional,  but the derivative
$$
\calA'_\kbf(\xbf)=(D_{x_j} \calA_\kbf(\xbf))_{j\in\Z},\qquad D_{x_j} \calA_\kbf(\xbf)= D_{x_j}(L_{k_{j-1}}(x_{j-1},x_{j})+L_{k_j}(x_j,x_{j+1}))
$$
is well defined.
A trajectory of the DLS is a pair $(\kbf,\xbf)\in K^\Z\times N^\Z$
such that $A'_\kbf(\xbf)=0$.
We call the trajectory $(\kbf,\xbf)$ admissible if the corresponding collision chain satisfies the jump condition (\ref{eq:jump}).

The Hessian
$$
\calA''_\kbf(\xbf)=(D_{x_i}D_{x_j} \calA_\kbf(\xbf))_{i,j\in\Z}
$$
of the action functional is 3-diagonal:
\begin{equation}
\label{eq:var}
\calA_\kbf''(\xbf)\ubf= \vbf,\qquad v_i=B_{i-1}u_{i-1} + A_i u_i + B_i^* u_{i+1},
\end{equation}
where
$$
B_i=D_{x_i}D_{x_{i+1}}L_{k_i}(x_i,x_{i+1})
$$
is the twist of the collision orbit $\gamma_i$.
The variational equation of the trajectory $(\kbf,\xbf)$ is $A_\kbf''(\xbf)\ubf=0$.
Under the twist condition, $B_i$ is invertible, and the variational equation defines the linear
Poincar\'e map $P_i:(u_{i-1},u_i)\to (u_i,u_{i+1})$.

For $n$-periodic collision chains, $\xbf$ is a critical point of the periodic action functional
\begin{equation}
\label{eq:per}
\calA_\kbf^{(n)}(\xbf)=\sum_{j=0}^{n-1}
L_{k_j}(x_j,x_{j+1}),\qquad \xbf=(x_1,\dots ,x_n),\quad
x_{n}=x_0 .
\end{equation}
We call the periodic collision chain nondegenerate if $\xbf$ is a nondegenerate  critical point of  $\calA^{(n)}_\kbf$.
If the twist condition holds, this is equivalent to the usual nondegeneracy condition  $\det (P-I)\ne 0$, where $P=P_n\circ\dots\circ P_1$ is  the linear  monodromy map.

Finite collision chains  joining the points $a,b\in M$ correspond to critical points of a finite sum
\begin{equation}
\label{eq:finite}
\calA_\kbf^{a,b}(\xbf)=\sum_{j=0}^{n}  L_{k_j}(x_j,x_{j+1}),\qquad x_0=a,\quad x_{n+1}=b,\quad \xbf=(x_1,\dots,x_n).
\end{equation}
We call the finite collision chain nondegenerate if the critical point $\xbf$ is nondegenerate.

Dynamics of the DLS is represented by the translation
$$
\calT:K^\Z\times N^\Z\to K^\Z\times N^\Z, \qquad (k_j,x_j)\to (k_{j+1},x_{j+1}).
$$
If $\calT$ has  a compact\footnote{The topology on $K^\Z\times N^\Z$ is the product topology.} 
invariant set $\Lambda\subset K^\Z\times N^\Z$ of  trajectories of the DLS, and the
collision map $\calF$ is well defined, then it will have a compact invariant set $\tilde \Lambda \subset K^\Z\times N^\Z$
with $\calF:\tilde\Lambda \to\tilde\Lambda $  topologically conjugate to $\calT:\Lambda\to\Lambda$.

The usual definition of a hyperbolic set is formulated in terms of the dichotomy of solutions
of the variational equation. It works under the twist condition, when the linear Poincar\'e maps $P_i$ are well defined.

A trajectory $(\kbf,\xbf)$ of the DLS is hyperbolic if for any  $j\in\Z$ there are stable and unstable subspaces $E_j^\pm\subset T_{x_{j}}N\times T_{x_{j+1}}N$  such that
$E_j^+\cap E_j^-=\{0\}$ and for any solution $\ubf=(u_i)$ of the variational equation
$w_j=(u_{j},u_{j+1})\in E^+_j$
 implies $w_i=(u_{i},u_{i+1})\in E^+_i$ for all $i>j$. Moreover $w_i$ decreases exponentially as $i\to\infty$: there is $C>0$ and $\lambda\in (0,1)$ such that
 $$
 \|w_i\|\le C\lambda^{i-j}\|w_j\|,\qquad i>j.
 $$
 Similarly for the unstable subspace:  $w_j=(u_{j},u_{j+1})\in E^-_j$
 implies $w_i=(u_{i},u_{i+1})\in E^-_i$ for all $i<j$ and $u_i$ decreases exponentially as  $i\to-\infty$:
 $$
 \|w_i\|\le C\lambda^{j-i}\|w_j\|,\qquad i<j.
 $$
 A compact $\calT$-invariant set $\Lambda$ of trajectories is hyperbolic if this holds for
 every trajectory $(\kbf,\xbf)\in\Lambda$ with $C,\lambda$ independent of the trajectory.

For our purposes another definition, not requiring the twist condition,  is more convenient.
If we use the Riemannian metric  to identify $T_{x_i}N$ and $T_{x_i}^*N$,
the Hessian $\calA''_\kbf(\xbf)$ becomes a linear operator
$\calA''_\kbf(\xbf):l_\infty \to l_\infty$, where $l_\infty$ is the Banach space of sequences
$$
\ubf=(u_i)_{i\in \Z},\qquad u_i\in T_{x_i}N,\quad \|\ubf\|_\infty=\sup_i\|u_i\|<\infty.
$$
If $\Lambda$ is a compact invariant set of the DLS, then  the Hessian
is a bounded operator: $\|\calA''_\kbf(\xbf)\|_\infty\le c=c(\Lambda)$ for any $(\kbf,\xbf)\in\Lambda$.

\begin{dfn}\label{def2}
We say that the trajectory $(\kbf,\xbf)$ is hyperbolic if
the Hessian  $\calA''_\kbf(\xbf)$  has bounded inverse in the $l_\infty$ norm.
We say that a compact $\calT$-invariant  set $\Lambda\subset K^\Z\times N^\Z$ of trajectories of the DLS is hyperbolic if this is true for all trajectories:
$\|\calA''_\kbf(\xbf)^{-1}\|_\infty\le C$ with $C=C(\Lambda)$ independent of the trajectory $(\kbf,\xbf)\in\Lambda$.
\end{dfn}

If the twist condition holds,   then,
as shown in \cite{Aubry-MacKay}, this definition of hyperbolicity
is equivalent to the standard one.\footnote{In  \cite{Aubry-MacKay} a single valued discrete Lagrangian was
considered, but in general the proof is the same.} But Definition \ref{def2} makes sense also without the twist condition,
for example when $N$
has connected components of different dimension, so the twist condition evidently fails.

\subsection{Main results}

Consider the system $(M\setminus N,H_\mu)$ with  Newtonian singularity on $N$ and the corresponding
degenerate billiard $(M,N,H_0)$ with Hamiltonian (\ref{eq:H0}).
Fix energy $E$.

\begin{thm}\label{thm:per}
Let   $\gamma$ be a nondegenerate periodic  collision chain    of the degenerate billiard
$(M,N,H_0=E)$.
There exists $\mu_0>0$ such that for any $\mu\in I_{\mu_0}=(-\mu_0,0)\cup(0,\mu_0)$ the    chain $\gamma$
is shadowed by  a periodic orbit $\gamma_\mu$ of the system $(M\setminus N,H_\mu=E)$.
\end{thm}

The shadowing error is of order $O(\mu\ln|\mu|)$, i.e.\ $d(\gamma_\mu(t),\gamma)\le c|\mu\ln|\mu||$.
At each near collision, the shadowing orbit $\gamma_\mu$ passes at a distance $\le c\mu$  from $N$.
However, for $\mu>0$ (attracting  singularity) it may have  collisions with $N$.
The regularized flow on the level $\{H_\mu=E\}$ has no singularity, so dynamics is always well defined.
If, for physical reasons, we need to avoid  regularizable  collisions, we have to impose an extra
condition on the  collision chain $\gamma=(\gamma_j)$.
Let
$$
v_j^\pm=\dot \gamma(t_j\pm 0)
$$
be the collision velocities at $j$-th collision point  $x_j=\gamma(t_j)$, and
let $u_j^\pm$ be their projections to the quotient space $T_{x_j}M/T_{x_j}N$.
The jump condition $\Delta p_j(t_j)\ne 0$ implies $u_j^+\ne u_j^-$.
For $\mu>0$ we assume  the no straight reflection condition $u_j^+\ne -u_j^-$:
\begin{equation}
\label{eq:straight}
v_j^+ + v_j^-\notin T_{x_j}N\quad \mbox{for all} \; j.
\end{equation}

Then the shadowing trajectory $\gamma_\mu$ will have no collisions: it passes   $N$ at the minimal distance
\begin{equation}
\label{eq:c12}
c_1\mu\le d(\gamma_\mu,N) \le c_2\mu,\qquad 0<c_{1}<c_2.
\end{equation}
Condition (\ref{eq:straight}) is less essential than  the jump condition (\ref{eq:jump})
since  dynamics is  well defined
also for trajectories colliding with $N$. For $\mu<0$  the  no straight reflection  condition is not needed.

\begin{rem}
If the twist condition holds (in particular all components of $N$ have the same codimension $d$),
then the periodic orbit $\gamma_\mu$ has $2d$  large Lyapunov exponents of order $O(\ln|\mu|)$.
Thus  $\gamma_\mu$ is strongly unstable, even if the corresponding periodic orbit of the DLS   is Lyapunov stable.
\end{rem}

Theorem \ref{thm:per} is a generalization of  a theorem in \cite{Bol-Neg:DCDS,Bol-Neg:RCD},
where it was proved for the case of second species solutions of the plane 3 body problem.

A similar statement holds for collision chains joining given points $a,b\in M\setminus N$.

\begin{thm}\label{thm:join}
Let   $\gamma$ be a nondegenerate   collision chain   of the degenerate billiard
$(M,N,H_0=E)$ joining the points $a,b\in M\setminus N$.
There exists $\mu_0>0$ such that for any $\mu\in I_{\mu_0}$ the    chain $\gamma$
is $O(\mu\ln|\mu|)$-shadowed by  an   orbit   of the system $(M\setminus N,H_\mu=E)$ joining $a,b$.
\end{thm}

The next theorem gives a hyperbolic invariant set of shadowing trajectories.

\begin{thm}
\label{thm:hyp}
Let $\Lambda\subset K^\Z\times N^\Z$ be a compact hyperbolic invariant set of the DLS
such that all orbits in $\Lambda$ are admissible.
There exists $\mu_0>0$ such that for any $\mu\in I_{\mu_0}$
and  any orbit $(\kbf,\xbf)\in\Lambda$
there exists a trajectory $\gamma_\mu$ of system $(M\setminus N,H_\mu=E)$  shadowing (as a non-parametrized curve)
the corresponding collision chain $\gamma$ of the degenerate billiard $(M,N,H_0=E)$.
Shadowing trajectories form a compact hyperbolic invariant set $\Lambda_\mu\subset \{H_\mu=E\}$ of system
$(M\setminus N,H_\mu=E)$.
\end{thm}

The shadowing error is of the same order  $O(\mu\ln|\mu|)$ as in Theorem \ref{thm:per}.
Recall that a trajectory $(\kbf,\xbf)$ of the DLS is admissible if the corresponding
collision chain satisfies the  jump condition  (\ref{eq:jump}).
For $\mu>0$ to avoid collisions we have to assume also the no straight reflections condition (\ref{eq:straight})
for trajectories in $\Lambda$. Then the shadowing trajectories satisfy (\ref{eq:c12}).

Note that in  Theorem \ref{thm:per} the periodic orbit of the degenerate billiard
does not need to be hyperbolic, so Theorems \ref{thm:per} and \ref{thm:hyp}
are formally independent.

To be honest, one of the main ingredients of the proof of Theorems \ref{thm:per}--\ref{thm:hyp}, Theorem \ref{thm:connect2},
will be proved only for $d=\codim N\le 3$.  The proof is based on Theorem \ref{thm:Shil} (the generalized  Shilnikov lemma) 
which holds for any codimension. However, to apply  Theorem \ref{thm:Shil}, we first need  to  regularize singularities.
We use the Levi-Civita regularization for $d\le 2$
and KS regularization \cite{K-S} for $d=3$. Collisions with $N$ (they are double collisions) are regularizable in any dimension,
but  standard multidimensional methods of regularization (e.g.\ Moser's regularization) are less convenient for our purposes
since regularization is not well defined in the limit $\mu\to 0$.
However, there is no doubt that   Theorem \ref{thm:connect2} is true for any $d$,
just the method of the proof needs to be changed. A multidimensional analog of the KS regularization is
the Clifford algebra  regularization which should give the proof of Theorem \ref{thm:connect2} for all $d>3$.
We  do not consider the case  $d>3$ since it has no applications in celestial mechanics
(unless one plans to do celestial mechanics in a space of dimension $>3$).

For a discrete scatterer $N$, Theorems  \ref{thm:per} and \ref{thm:hyp} were proved in \cite{Bol-Mac:centers}
and used to prove the existence of chaotic second species solutions of the restricted circular  3 body problem.
A version of these theorems for the elliptic restricted 3 body problem was proved in \cite{Bol:DCDS}
(then $N$ is one-dimensional). A version of Theorem  \ref{thm:per} was proved in
\cite{Bol-Neg:DCDS} for the plane nonrestricted 3 body problem.
Then $N$ is 2-dimensional but becomes 1-dimensional after reduction of symmetry.

\subsection{Shadowing for systems with symmetry}

\label{sec:reduce}

Formally   Theorems \ref{thm:per} and \ref{thm:hyp} are of little use  in  celestial mechanics.
Indeed,   Hamiltonian systems of celestial mechanics usually have translational or rotational symmetry and so
they do not possess nondegenerate
periodic orbits or hyperbolic invariant sets. Hence Theorems \ref{thm:per} and \ref{thm:hyp}  do not apply.
The exception is Theorem \ref{thm:join}: it works  also in the presence of symmetry.
Indeed, symmetry is broken by fixing the end points of a trajectory (if they are not fixed points of the group action),
so nondegenerate connecting chains may exist. Restricted    problems of celestial mechanics also have symmetry broken
and then all Theorems \ref{thm:per}--\ref{thm:hyp} work.

To apply  Theorems \ref{thm:per} and \ref{thm:hyp} in celestial mechanics, we have to  reduce symmetry.
We describe the reduction in the  simplest situation arising in applications, see also \cite{Bol-Neg:DCDS}.
Suppose the degenerate billiard $(M,N,H)$ has an abelian  symmetry group $\A^s$,
where $\A^s$ is a torus $\T^s=\R^s/\Z^s$, or $\R^s$,
or their product (cylinder).
 More precisely, suppose there is a smooth group action  $\Phi_\theta:M\to M$,
$\theta\in\A^s$, which preserves the Hamiltonian and the scatterer:
$$
\Phi_\theta (N)=N,\qquad H(\Phi_\theta (q),p)=H(q,D\Phi_\theta(q)^*p).
$$

For any $\xi\in\R^s$, the one-parameter
symmetry group $\Phi_{t\xi}$ is generated by the vector field  $u_\xi(q)=X(q)\xi$, where
$$
X(q)=D_\theta\big|_{\theta=0}\Phi_\theta(q):\R^s\to T_qM.
$$
Let
$$
G_\xi(q,p)=\langle u_\xi(q),p\rangle
$$
be the corresponding Noether  integral \cite{Arnold,AKN} of the Hamiltonian system.
Then
$$
G:T^*M\to (\R^s)^*,\qquad \langle G(q,p),\xi\rangle=G_\xi(q,p),\qquad \xi\in\R^s,
$$
is the  momentum integral.
Since $u_\xi$ is tangent to $N$, $G$ is preserved by the reflection and so it will
be also an integral of the degenerate billiard $(M,N,H)$.

The corresponding  DLS with the Lagrangian $\calL=\{ L_k\}_{k\in K}$ has the symmetry
$$
 L_k(\Phi_\theta (x_-),\Phi_\theta (x_+))= L_k(x_-,x_+).
$$
The action functional (\ref{eq:Ak}) is   invariant:
$$
\calA_\kbf(\Phi_\theta \xbf)=\calA_\kbf(\xbf).
$$
Thus for any $\xi\in\R^s$, $\ubf_\xi=(u_\xi(x_j))_{j\in\Z}$ is  in the kernel of the Hessian
$\calA''_\kbf(\xbf)$,
and the Hessian is non-invertible: there are no nondegenerate periodic orbits
or hyperbolic trajectories except fixed points of the group $\Phi_\theta$.

We call an $n$-periodic collision chain $\gamma=(\gamma_i)_{i\in\Z}$
nondegenerate modulo symmetry if it has only degeneracy coming from symmetry.  
The corresponding critical point  $\xbf$
of the action functional (\ref{eq:per}) satisfies
$$
D^2\calA^{(n)}_\kbf(\xbf)\vbf=0\quad\Rightarrow\quad v_i=u_\xi(x_i),\qquad \xi\in \R^s.
$$

Suppose now that the system $(M\setminus N,H_\mu)$  with Newtonian singularities has a symmetry group:
$$
 H_\mu(\Phi_\theta (q),p)=H_\mu(q,D\Phi_\theta(q)^*p),\qquad \theta\in\A^s.
$$
Then $\Phi_\theta$ is a symmetry group   of the corresponding degenerate billiard $(M,N,H_0)$.
 We have the following version of Theorem \ref{thm:per} for systems with symmetry.

\begin{thm}\label{thm:symm}
Let   $\gamma$ be a nondegenerate modulo symmetry periodic collision  chain
of the degenerate billiard $(M,N,H_0=E)$.
There exists $\mu_0>0$ such that for any $\mu\in  I_{\mu_0}$ the    chain $\gamma$
is shadowed by  a periodic orbit $\gamma_\mu$ of the system $(M\setminus N,H_\mu=E)$.
\end{thm}

Of course $\gamma_\mu$ is defined modulo symmetry $\gamma_\mu\to\Phi_\theta\gamma_\mu$.
In Theorem \ref{thm:symm} it is not possible to prescribe the value of the momentum integral $G$ of the periodic orbit $\gamma_\mu$.
To find trajectories with given value of $G$, we need to consider orbits periodic modulo symmetry:
$\gamma(t+T)=\Phi_\theta\gamma(t)$.

The   discrete action functional (\ref{eq:per}) is modified as follows:
$$
\calP_\kbf(\xbf,\theta)=\calA^{(n)}_\kbf(\xbf)-\langle G,\theta\rangle,\qquad \xbf=(x_1,\dots,x_n),\quad \theta\in\R^s,\quad x_n=\Phi_\theta (x_0).
$$
Critical points $(\xbf,\theta)$ of $\calP_\kbf$ correspond to  collision chains $\gamma$ which are periodic modulo symmetry and have integral $G$.
We call   $\gamma$ nondegenerate if $(\xbf,\theta)$ is a nondegenerate (modulo symmetry)
critical point of $\calP_\kbf$.

\begin{thm}
\label{thm:symm-G}
Let   $\gamma$ be a nondegenerate periodic modulo symmetry collision chain
of the degenerate billiard $(M,N,H=E)$. Let $G$ be its  momentum integral.
There exists $\mu_0>0$ such that for any $\mu\in  I_{\mu_0}$ the chain $\gamma$
is shadowed modulo symmetry by  a periodic modulo symmetry orbit $\gamma_\mu$ of the system $(M\setminus N,H_\mu=E)$
with the momentum integral $G$.
\end{thm}

Shadowing modulo symmetry means that $d(\Phi_{\theta(t)}\gamma_\mu(t),\gamma)\le c|\mu\ln |\mu||$  for some $\theta(t)\in\A^s$.

To prove Theorem \ref{thm:symm-G},
 we perform symmetry reduction.
 Suppose that the quotient space
$\tilde M=M/\Phi_\theta$ is a smooth manifold and the projection $\pi:M\to\tilde M$
is a smooth fiber bundle with fiber $\A^s$.
For simplicity   assume that the fibre bundle $\pi:M\to\tilde M$ is trivial. This is always true locally.
Then $\tilde M$ can be realized as a cross section $\tilde M\subset M$ of the group action $\Phi_\theta$.

Let $L$ be the Lagrangian  (\ref{eq:L}).
Define the reduced Lagrangian (Routh function) on $T\tilde M$ by
\begin{equation}
\label{eq:Routh}
\tilde L(q,\dot q)=\Crit_\xi(L(q,\dot q+u_\xi(q))-\langle G,\xi\rangle),\qquad q\in \tilde M,\quad \dot q\in T_{q}\tilde M,
\end{equation}
where $\Crit_\xi$ means taking a critical value with respect to $\xi\in\R^s$.
Since $L$ is convex in the velocity, the Routh function is well
defined. For  the  standard  definition see \cite{Arnold,AKN}.

Let $\tilde H$ be the Hamiltonian corresponding to $\tilde L$.
Then trajectories of the Hamiltonian system $(M,H)$ with the momentum $G$
are projected to trajectories of the reduced  Hamiltonian system $(\tilde M,\tilde H)$.

If $(M,N,H)$ is a degenerate billiard with symmetry, then the reduced degenerate billiard is
$(\tilde M,\tilde N,\tilde H)$, where $\tilde N=N/\Phi_\theta$.

If the system with singularities $(M\setminus N,H_\mu)$ has a symmetry $\Phi_\theta$, then for fixed momentum $G$,
the degenerate billiard corresponding to
the reduced system $(\tilde M\setminus \tilde N,\tilde H_\mu)$ will  be  the reduced billiard $(\tilde M,\tilde N,\tilde H_0)$.
Now we can apply Theorem \ref{thm:per} to the reduced system with singularities and  to the corresponding reduced billiard.
This proves Theorem \ref{thm:symm-G}.

\begin{rem}
If the group $\A^s$ is a torus  or a cylinder, then the fibration is nontrivial in general.
Then the construction of the  reduced   system is local.
The global version needs choosing a connection for the fibre bundle $\pi:M\to\tilde M$
and using the symplectic structure  twisted by the curvature form of the connection \cite{AKN}.
However this is not needed  in this paper since our results are essentially local:
we can  assume that the collision chains
lie in a domain $U\subset M$ such that  fibre bundle $\pi:U\to \tilde U$ is trivial.
\end{rem}

We can perform  symmetry reduction also for the DLS describing the billiard.
For any trajectory $(\kbf,\xbf)$,
the Noether integral corresponding to $\xi\in\R^s$ is
$$
G_\xi=\langle u_\xi(x_j),y_j\rangle,
$$
where $y_j$ is the momentum (\ref{eq:y_j}).

The fibration $\pi:M\to \tilde M$ defines a fibration $\pi:N\to \tilde N$ to the orbits of the group action.
We assume that it is trivial. Then $\tilde N=N/\Phi_\theta$ can be identified with a  cross section $\tilde N\subset N$ of the group action $\Phi_\theta|_N$.

For a fixed value of the integral   $G$, define the reduced discrete Lagrangian (discrete Routh function) by the Legendre transform
\begin{equation}
\label{eq:Routh}
\tilde  L_k(x_-,x_+)=
\Crit_\theta(L_k(x_-,\Phi_\theta(x_+))-\langle G,\theta\rangle),\qquad x_\pm\in \tilde N .
\end{equation}
  This requires a  twist condition: the bilinear form
\begin{equation}
\label{eq:nondeg}
\xi,\eta\to \langle B_k(x_-,x_+)u_\xi(x_-),u_\eta(x_+)\rangle,\qquad \xi,\eta\in\R^s,
\end{equation}
is nondegenerate. Here  $B_k(x_-,x_+)=D_{x_-}D_{x_+}L_k(x_-,x_+)$ is the twist of the Lagrangian $L_k$.
In general the reduced discrete Lagrangian  is locally defined: it is a function on an open set
$\tilde U_k\subset \tilde N\times \tilde N$.

For any trajectory $(\kbf,\xbf)$ of the DLS with momentum integral $G$ setting  $\tilde \xbf=(\tilde x_j)$, $\tilde x_j=\pi(x_j)$,
we obtain a trajectory $(\kbf,\tilde\xbf)$ of the reduced DLS with the Lagrangian $\tilde\calL=\{\tilde L_k\}_{k\in K}$.
Conversely, a trajectory $(\kbf,\tilde\xbf)$ of the reduced DLS defines a (nonunique) trajectory $(\kbf,\xbf)$ of the original DLS
with  momentum $G$.

Now we can apply Theorem \ref{thm:hyp}  to the reduced Hamiltonian system $(\tilde M\setminus\tilde N,\tilde H_\mu)$ and the corresponding
reduced DLS   and obtain the existence of hyperbolic modulo symmetry
invariant  sets on a level set of $G$  for the system $(M\setminus N,H_\mu=E)$ when the corresponding reduced  DLS has
a compact hyperbolic invariant set.

In the next publication these results will be used to study chaotic second species solutions
of the nonrestricted 3 body problem.

\section{Proofs}

\label{sec:proof}

In this section we prove   Theorems \ref{thm:per} and \ref{thm:hyp}.
The proof of Theorem \ref{thm:join} is similar. The proofs are based on a local connection result --
Theorem \ref{thm:connect} which is proved in  section \ref{sec:reg}.

\subsection{Local connection}

\label{sec:local}

Let $d$ be the distance in $M$ defined by the Riemannian metric $\|\,\cdot\,\|=\|\,\cdot\,\|_0$.
We parameterize a tubular neighborhood
$$
N_\rho=\{q\in M:d(q,N)\le \rho\}
$$
by the exponential map
$$
f:T^\perp N\to M,\qquad q=f(x,u)=\exp_xu,\qquad x\in N,\; u\in T_x^\perp N.
$$
Let $D\Subset N$ be an  open set with  compact closure.
Then for small $\rho>0$,
$$
U_\rho =\{q=f(x,u):x\in D,\; \|u\|\le \rho\}\subset N_\rho
$$
has  smooth boundary
\begin{equation}
\label{eq:Sigma}
\Sigma_\rho =\{q=f(x,u):x\in D,\; \|u\|= \rho\}
\end{equation}
 and $f$  is a diffeomorphism onto  $U_\rho$. For $q\in U_\rho$ we have $d(q,N)=\|u\|$.

Consider a  degenerate billiard $(M,N,H_0=E)$. Suppose that $D\Subset N\cap \calD_E$ is
contained in the domain of possible motion (\ref{eq:DE}).
There exists $r>0$ such that for any $x_0\in D$ we have $B_r(x_0)\subset \calD_E$
and for any pair of points $q_-\ne q_+$ in the ball $B_r(x_0)$ there exists
a trajectory  $\gamma$  of system $(M,H=E)$ (geodesic of the Jacobi metric) joining $q_\pm$
in $B_r(x_0)$. The trajectory $\gamma$ smoothly depends on $q_+\ne q_-$.
Let $S(q_-,q_+)=J_0(\gamma)$ be its Maupertuis action (\ref{eq:S}).

Fix arbitrary large\footnote{We denote by $c,C$ several large fixed constants.}
$C>0$ and  let
\begin{equation}
\label{eq:Prho}
P_\rho=\{(q_+,q_-)\in\Sigma_\rho^2: d(q_+,q_-)\le C\rho\}.
\end{equation}
We will  connect a pair of points $q_+,q_-\in P_\rho$ by a  billiard trajectory of energy $E$
having a single reflection from $N$ at a point $x_0$.

\begin{prop}
\label{prop:reflect}
Let $\rho=\rho(C,D)>0$ be sufficiently small.
Then for any $(q_-,q_+)\in P_\rho$:
\begin{itemize}
\item
There exists $x_0\in N$  and a  trajectory $\gamma$ of the degenerate billiard $(M,N,H=E)$
joining $q_+,q_-$ in $B_r(x_0)\cap N_\rho$ after a  reflection at  $x_0$.
Thus $\gamma=\gamma_+\cdot \gamma_-$ is a concatenation
of a trajectory $\gamma_+$ which joins
$q_+$ with $x_0$ and $\gamma_-$ which joins $x_0$ with $q_-$.
\item $x_0=\xi(q_+,q_-)$ and $\gamma$ smoothly depend on $(q_+,q_-)\in P_\rho$.
\item The Maupertuis action $J_0(\gamma)=R_0(q_+,q_-)$ is a smooth function on $P_\rho$ and
\begin{equation}
\label{eq:g0}
R_0(q_+,q_-)=\Crit_{x\in N}R(q_+,x,q_-),\qquad R=S(q_+,x)+S(x,q_-).
\end{equation}
More precisely, $x_0$ is the only critical point of $x\to R(q_+,x,q_-)$ in $N\cap B_r(x_0)$, and it is nondegenerate.
\end{itemize}
\end{prop}

Let $p_\pm\in T_{q_\pm}^*\Sigma_\rho$ be the momenta of $\gamma$ at $q_\pm\in\Sigma_\rho$.
Then $R(q_+,x,q_-)$ is the generating function of the Lagrangian relation $\calR$ between the points
$(q_+,p_+)$ and $(q_-,p_-)$. Note that $\calR$ is not a  map unless $N$ is a hypersurface (ordinary billiard), then
$\calR:(q_+,p_+)\to (q_-,p_-)$ is a symplectic map of a set in $T^*\Sigma_\rho$.

Proposition \ref{prop:reflect}  is a  familiar  property of systems with elastic reflections (a version of Fermat' principle).
However we  give a proof since the notations will be needed in the next theorem.

Since Proposition \ref{prop:reflect} is local:
all trajectories lie in a neighborhood of some point $x_0\in N$,
without loss of generality we may assume that $D\subset N$ is contractible and is contained in a coordinate chart in $N$.
Then the normal bundle $T^\perp N$ is trivial over $D$ and
we can  choose an orthonormal basis $e_1(x),\dots,e_d(x)$ in $T_x^\perp N$ smoothly depending on $x\in D$.
Then the exponential map
\begin{equation}
\label{eq:coord}
q=f(x,u),\quad u=u_1e_1+\dots+u_de_d,
\end{equation}
defines coordinates $x\in D$, $u\in B_\rho=\{ u\in \R^d:|u|<\rho\}$, in $U_\rho$. 
Then $q\in\Sigma_\rho$ when $u\in S_\rho=\partial B_\rho$.
We denote by $v\in\R^d$ the momentum conjugate to $u$ and by $y$ the momentum conjugate to $x\in D$.
The a trajectory of the Hamiltonian system is represented by $z(t)=(x(t),y(t))$, $u(t)$, $v(t)$.

Let $F_0$ be the Hamiltonian  (\ref{eq:F}) corresponding to the Hamiltonian $H_0$:
\begin{equation}
F_0(x,y)=\frac12\langle A_0(x)(y-a_0(x)),y-a_0(x)\rangle+W_0(x),
\label{eq:F}
\end{equation}
and let
\begin{equation}
\label{eq:calK}
\calK=\{(x,y)\in T^*D: F_0(x,y)\le E-\eps\}\subset \calM_E.
\end{equation}
If $\rho>0$ is small enough, for any $z_0=(x_0,y_0)\in\calK$ and any $u_-\in S_\rho=\partial B_\rho$  there is   $t_->0$  and a trajectory
$\gamma_-=\gamma_-(z_0,u_-):[0,t_-]\to U_\rho$ with $H=E$  satisfying the boundary conditions
 \begin{equation}
\label{eq:gamma-}
z(0)=z_0,\quad u(0)=0,\quad u(t_-)=u_-.
\end{equation}
Similarly, for any $u_+\in S_\rho$ there  is $t_+<0$ and a trajectory $\gamma_+(z_0,u_+):[t_+,0]\to U_\rho$
with $H=E$ satisfying the boundary conditions
 \begin{equation}
\label{eq:gamma+}
z(0)=z_0,\quad u(0)=0,\quad u(t_+)=u_+ .
\end{equation}
The concatenation $\gamma_+\cdot \gamma_-$ is a reflection trajectory of the degenerate billiard with collision point $x_0$
and tangent collision momentum $y_0$.

Indeed, for $u=0$   we have
$$
H_0(x,y,0,v)=F_0(x,y)+\frac12 |\dot u|^2=E.
$$
For a solution of the  Hamiltonian system
with the initial condition
$z(0)=z_0$, $u(0)=0$ and energy $E$, we have
$$
 u(t)=u(z_0,\dot u(0),t)=t\dot u(0)+O(t^2),\qquad |\dot u(0)|=\sqrt{2(E-F_0(z_0))}=\nu(z_0).
 $$
For small $\rho>0$ the equation $u(t_\pm)=u_\pm$  can be solved for
 \begin{eqnarray*}
 t_\pm&=&t_\pm(z_0,u_\pm)=\mp \frac{\rho}{\nu(z_0)}+O(\rho^2),\\
 \dot u(0)&=&\dot u(z_0,u_\pm)=\mp  \nu(z_0)e_\pm+O(\rho),\qquad  u_\pm=\rho e_\pm,\quad |e_\pm|=1.
 \end{eqnarray*}
Here $O(\rho)$ means a function of the form $\rho h(z_0,e_\pm,\rho)$ where $h$ is $C^1$ bounded  as $\rho\to 0$.
The corresponding trajectories $\gamma_\pm$ satisfy (\ref{eq:gamma+})--(\ref{eq:gamma-}).

In local coordinates in $D$,  we have $\gamma_\pm(t_\pm)=f(x_\pm,u_\pm)$, where
\begin{equation}
\label{eq:dotx}
x_\pm=\xi_\pm(z_0,u_\pm)=x_0+t_\pm \dot x(0)+O(\rho^2),\qquad \dot x(0)=A_0(x_0)(y_0-a(x_0)).
\end{equation}
To prove Proposition \ref{prop:reflect}, for given $q_\pm=f(x_\pm,u_\pm)$ such that $(q_+,q_-)\in P_\rho$,
we need to find $z_0=(x_0,y_0)$  such that
\begin{equation}
\label{eq:xi}
\xi_\pm(z_0,u_\pm)=x_\pm.
\end{equation}
We have $|x_+-x_-|\le c\rho$ with $c>0$ independent of $\rho$.
Using (\ref{eq:dotx}), equations (\ref{eq:xi})  can be rewritten as
$$
x_0=\frac 12(x_++x_-)+O(\rho^2),\quad  y_0=a(x_0)-\frac{\nu(z_0)}{2\rho} A_0^{-1}(x_0)(x_+-x_-)+O(\rho).
$$
For small $\rho>0$,  equations (\ref{eq:xi}) satisfy the condition of the implicit function theorem
and so they can be solved for $(x_0,y_0)=z_0(q_+,q_-)$.

Proposition \ref{prop:reflect} is proved.
\qed

\medskip

Next we  formulate a similar local connection result for the system $(M\setminus N,H_\mu=E)$ with Newtonian singularities.
The connection trajectory will be close to the reflection trajectory $\gamma_+\cdot \gamma_-$ of the degenerate billiard $(M,N,H_0=E)$  
in Proposition \ref{prop:reflect}. We need another restriction on the points $q_\pm$ we try to connect.
We write it in local coordinates defined in (\ref{eq:coord}).

Fix small $\delta>0$ and let
\begin{equation}
\label{eq:Qrho}
Q_\rho=\{(q_+,q_-)\in P_\rho: q_\pm=f(x_\pm,u_\pm),\; |u_++u_-|\ge \delta\rho\}.
\end{equation}
Thus we do not want the points $q_\pm$ to  be nearly opposite with respect to $N$.

\begin{thm}\label{thm:connect}
Let $\rho=\rho(\delta,C,D)>0$ be sufficiently small.
There exists $\mu_0>0$ such that for all $(q_+,q_-)\in Q_\rho$ and $\mu\in I_{\mu_0}= (-\mu_0,0)\cup (0,\mu_0)$:
\begin{itemize}
\item
There exists  a unique (up to a time shift) trajectory $\alpha_\mu$ of system $(M\setminus N,H_\mu=E)$ joining $q_+$ and $q_-$ in $ B_r(x_0)$, where $x_0=\xi(q_+,q_-)$.
\item
$\alpha_\mu$ smoothly depends on $(q_+,q_-,\mu)\in Q_\rho\times I_{\mu_0}$
and uniformly converges (as a nonparametrized curve) as $\mu\to 0$ to the billiard trajectory  $\gamma_+\cdot \gamma_-$ in Proposition
\ref{prop:reflect}.
\item
The minimal distance $d(\alpha_\mu,N)$  is attained at a point $q_\mu=f(x_\mu,u_\mu)$,
 which converges to $x_0=\xi(q_+,q_-)$ as $\mu\to 0$:
$$
|u_\mu|\le c|\mu|,\quad d(x_\mu,x_0)\le c|\mu\ln|\mu||.
$$
\item
The Maupertuis action of $\alpha_\mu$  has the form
\begin{equation}\label{eq:Jmu}
J_{\mu}(\alpha_\mu)=\int_{\alpha_\mu}p\, dq
=R_\mu(q_+,q_-)=R_0(q_+,q_-)+O(\mu\ln|\mu|),
\label{eq:gmu}
\end{equation}
where $R_0(q_+,q_-)$ is the  action (\ref{eq:g0}) of the billiard trajectory $\gamma_+\cdot\gamma_-$ and $O(\mu\ln|\mu|)$
means a function $h$ such that
$$
\|h\|_{C^2(Q_\rho)}\le c |\mu\ln|\mu||
$$
with a  constant $c$ independent of $\mu$.
 \end{itemize}
\end{thm}

Theorem \ref{thm:connect} implies that the  symplectic  map $\calP_\mu:(q_+,p_+)\to (q_-,p_-)$ of $T^*\Sigma_\rho$,
which has no limit as $\mu\to 0$, does have a smooth limit  if represented as a
Lagrangian relation with the generating function $R_\mu$.

\begin{rem}
For the attracting force ($\mu>0$) the connecting trajectory $\alpha_\mu$ in Theorem \ref{thm:connect} may have a
regularizable collision with $N$ (although the set of $(q_+,q_-)$
with this property is negligible).  To avoid this,  we have to replace $Q_\rho$
with the  set
$$
\hat Q_\rho=\{(q_-,q_+)\in Q_\rho:q_\pm=f(x_\pm,u_\pm),\; |u_+-u_-|\ge \delta\rho
 \}.
$$
If $(q_+,q_-)\in \hat Q_\rho$, then the billiard trajectory $\gamma_+\cdot \gamma_-$ satisfies the no straight reflection
condition (\ref{eq:straight}) at $x_0$.
Then the shadowing orbit $\alpha_\mu$ will  satisfy
$$
c_1\mu\le d(\alpha_\mu,N)\le c_2\mu,\qquad c_{1,2}>0.
$$
\end{rem}

Theorem \ref{thm:connect} is a  generalization of the result proved in \cite{Bol-Neg:RCD}.
We will  deduce it from the following Theorem \ref{thm:connect2}. Let $\calK\subset \calM_E$ be  the set  (\ref{eq:calK}).

\begin{thm}\label{thm:connect2} Fix $\delta>0$.
Let $\rho>0$ be sufficiently small. There exists $\mu_0>0$ such that
for all $\mu\in I_{\mu_0}$, any $z_0=(x_0,y_0)\in \calK$ and any
$u_\pm\in S_\rho$ such that  $|u_++u_-|\ge \delta \rho$:
\begin{itemize}
\item
There exists  a trajectory $\gamma_\mu:[t_+,t_-]\to
N_\rho$, $t_+<0<t_-$, of system $(M\setminus N,H_\mu=E)$   satisfying the initial-boundary conditions
$u(t_\pm)=u_\pm$, $z(0)=z_0$.
\item
$\gamma_\mu$ smoothly depends on $(z_0,u_+,u_-,\mu)\in \calK\times
S_\rho^2\times I_{\mu_0}$ and converges, as $\mu\to 0$, to a
trajectory $\gamma_+\cdot\gamma_-$ of the degenerate billiard having  a reflection from $N$ at $x_0$
with the tangent momentum $y_0$.
\item
The Maupertuis action of $\gamma$ has the form
\begin{eqnarray}
\label{eq:Jmu}
J_\mu(\gamma_\mu)  =\psi_+(z_0,u_+)+\psi_-(z_0,u_-)+O(\mu\ln |\mu|),
\end{eqnarray}
where $\psi_\pm(z_0,u_\pm)$  are the Maupertuis actions of the trajectories $\gamma_\pm(z_0,u_\pm)$
of the Hamiltonian system $(M,H_0=E)$
satisfying the boundary conditions (\ref{eq:gamma+})--(\ref{eq:gamma-}).
\item The end points of $\gamma_\mu$ satisfy
\begin{eqnarray}
x(t_\pm)
=x_\mu^\pm(z_0,u_+,u_-)=\xi_\pm(z_0,u_\pm)+ O(\mu\ln |\mu|).
\label{eq:xpm}
\end{eqnarray}
Here $O(\mu\ln |\mu|)$ means a function which is    uniformly $C^1$ bounded on
$\calK\times  S_\rho^2$ for $\mu\in I_{\mu_0}$ by $c|\mu\ln |\mu||$.
\item
If $\mu<0$, or $\mu>0$ and   $|u_+-u_-|\ge \delta \rho$, then
\begin{equation}
\label{eq:min} \mu c_1\le d(\gamma_\mu,N)=\min|u(t)|\le
\mu c_2, \qquad 0<c_1<c_2.
\end{equation}
\end{itemize}
\end{thm}

Let us deduce Theorem \ref{thm:connect} from Theorem \ref{thm:connect2}.
We have $q_\pm=f(x_\pm,u_\pm)$, $u_\pm\in S_\rho$, where $d(x_+,x_-)\le c\rho$ and $|u_++u_-|\ge \delta\rho$.
We need to find $z_0\in\calK$ such that
the trajectory $\gamma_\mu$  in Theorem \ref{thm:connect2}
corresponding to $u_\pm$ and $z_0\in\calK$ satisfies
$x_\mu^\pm(z_0,u_-,u_+)=x_\pm$.

For $\mu=0$ this is done in the proof of Proposition
\ref{prop:reflect} and $z_0=z_0(q_+,q_-)$ was obtained as a nondegenerate solution of equations
(\ref{eq:xi}). Since the implicit function
theorem worked  for $\mu=0$, by (\ref{eq:xpm}),
for small $\mu$ it will work also here. \qed

\medskip

Theorem \ref{thm:connect2} is proved (for $d\le 3$) in  section \ref{sec:reg}.

\subsection{Proof of Theorem \ref{thm:per}}

 The idea  of the proof is to represent the shadowing trajectory  of system $(M\setminus N,H_\mu=E)$
as a critical point of a functional $\Phi_\mu$ which is nonsingular as $\mu\to 0$.

Let   $\gamma^0=(\gamma_j^0)$ be  a nondegenerate   $n$-periodic collision chain of the degenerate billiard $(M,N,H_0=E)$.
Suppose the collision orbit $\gamma_j^0$ connects
the points $x_j^0\in N$ and $x_{j+1}^0\in N$.
There is a $n$-periodic sequence $\kbf=(k_j)$ such that $\xbf^0=(x_j^0)$ is a nondegenerate critical point of the function (\ref{eq:per}).

Since $x_j^0\in\calD_E$, there exists $r>0$ such that $B_r(x_j^0)\Subset\calD_E$ for all $j$.
Set $D_j=B_r(x_j^0)\cap N$ and  $D=\cup D_j$.

Take small $\rho>0$. Suppose that the collision orbit $\gamma_j^0$ crosses $\Sigma_\rho$
at the points $s_j^-$ near $x_j^0$ and $s_{j+1}^+$ near $x_{j+1}^0$.
 Since   $\gamma_j^0$
is not tangent to $N$ at the end points, taking $\rho$ small enough we may assume that
there is a constant $C>0$, independent of $\rho$, such that
$$
d(s_j^+,s_j^-)<C\rho.
$$
Hence $(s_j^+,s_j^-)\in P_\rho$, where $P_\rho$ is the set (\ref{eq:Prho}) corresponding to $C,D$.
We take $\rho>0$ so small that Proposition \ref{prop:reflect} holds in $P_\rho$.
There is $\eps>0$ such that for
$$
q_j^\pm \in B_j^\pm=\Sigma_\rho\cap B_\eps(s_j^\pm),
$$
we have $(q_j^+,q_j^-)\in P_\rho$. Then the  action function $R_0$ in Proposition \ref{prop:reflect}
is defined on  $B_j^+\times B_j^-$.

In the   coordinates  $x\in D_j$, $u\in B_\rho$ in a neighborhood of $x_j^0$, we have $s_j^\pm=f(x_j^\pm,u_j^\pm)$.
 The jump condition implies that there is $\delta>0$ such that $|u_j^-+u_j^+|\ge 2\delta \rho$.
Then if $\eps>0$ is small enough, $q_j^\pm\in B_j^\pm$ implies $(q_j^+,q_j^-)\in Q_\rho$,
where $Q_\rho$ is the set (\ref{eq:Qrho}) corresponding to $C,D,\delta$.
By Theorem \ref{thm:connect}, if $\rho>0$ is small enough, there is $\mu_0>0$ such that
for $\mu\in I_{\mu_0}=(-\mu_0,0)\cup (0,\mu_0)$ there exists an orbit $\alpha_j=\alpha_j(q_j^+,q_j^-,\mu)$
of system $(M\setminus N,H_\mu=E)$
joining $q_j^+$ and $q_j^-$ in $N_\rho \cap B_r(x_j^0)$. Its   action $J_\mu(\alpha_j)=R_\mu(q_j^+,q_{j}^-)$ is given by (\ref{eq:gmu}).

Since the points $x_j^0$ and $x_{j+1}^0$ are non-conjugate along $\gamma_j^0$,  for small $\rho>0$
the points $s_j^-$ and $s_{j+1}^+$ are not conjugate along the corresponding segment of $\gamma_j^0$.
Hence  there exist $\eps>0$
and $\mu_0>0$
such that for any $\mu\in (-\mu_0,\mu_0)$, any points $q_j^-\in B_j^-$ and $q_{j+1}^+\in B_{j+1}^+$
are joined by a unique trajectory $\beta_j=\beta_j(q_j^-,q_j^+,\mu)$
of system $(M\setminus N,H_\mu=E)$ which is close to $\gamma_j^0$.
Let $F_j(q_j^-,q_{j+1}^+,\mu)=J_\mu(\beta_j)$ be its Maupertuis action.

Consider the function
\begin{equation}
\label{eq:Phi}
\Phi_\mu(\qbf)=\sum_{j=1}^n  (F_j(q_j^-,q_{j+1}^+,\mu)+R_\mu(q_j^+,q_{j}^-)),\qquad
q_{n+1}^\pm=q_1^\pm,
\end{equation}
where
\begin{equation}
\label{eq:q}
\qbf=(q_1^+,q_1^-,\dots, q_n^+,q_n^-)\in \calB=B_1^+\times B_1^-\times\dots\times B_n^+\times B_n^-.
\end{equation}
Then $\Phi_\mu(\qbf)$ is the Maupertuis action $J_\mu(\hat \gamma)$ of the concatenation $\hat\gamma$  of the trajectories $\alpha_j$, $\beta_j$
defined above. This is a broken trajectory with momentum discontinuous at $q_j^\pm$.

\begin{lem}
If $\qbf\in\calB$ is a critical point of $\Phi_\mu$,
then the concatenation  $\hat\gamma$ is a smooth
periodic trajectory of system $(M\setminus N,H_\mu=E)$.
\end{lem}

Indeed, by Hamilton's first variation formula,
$$
\delta\Phi_\mu(\qbf)=\delta J_\mu(\hat\gamma)=\sum_{j=1}^n(\langle \Delta p_j^+,\delta q_j^+\rangle+\langle \Delta p_j^-,\delta q_j^-\rangle)=0,\qquad \delta q_j^\pm\in T_{q_j^\pm}\Sigma_\rho,
$$
where $\Delta p_j^\pm$ is the jump of the momentum at $q_j^\pm$.
Hence $\Delta p_j^\pm\perp T_{q_j^\pm}\Sigma_\rho$. Since the Hamiltonian $H_\mu=E$ has no jump, and $\Sigma_\rho$ is a hypersurface,
this implies $\Delta p_j^\pm=0$.

Indeed, let $u_j$ be the initial  velocity of $\alpha_j$ at $q_j^+$ and $v_j$ the final velocity of $\beta_{j-1}$ at $q_j^+$. Then
$\Delta v_j=u_j-v_j $ is orthogonal to $T_{q_j^+}\Sigma_\rho$ with respect to the Riemannian metric and
$$
\|u_j\|^2=\|v_j\|^2=2(E-W_\mu(q_j^+)).
$$
This implies that either $u_j=v_j$ and  $\Delta v_j= 0$, so the concatenation $\beta_{j-1}\cdot \alpha_j$ is smooth at $q_j^+$,
or the concatenation has an elastic reflection from $\Sigma_\rho$,
and then $\Delta v_j\ne 0$. The second case is impossible since  $\alpha_j\subset N_\rho$,
so its velocity $u_j$ at $q_j^+$ points outside $\Sigma_\rho$,  and
the velocity $v_j$ of $\beta_j$ at $q_j^+$  is close to the velocity of $\gamma_{j-1}^0$ at $s_j^+$
so it also points outside $\Sigma_\rho$.

Hence   the concatenation $\hat\gamma$ is smooth at $q_j^+$. Similarly for $q_j^-$.
\qed

\medskip

Let us show that for $\mu=0$ the function $\Phi_0$ has a nondegenerate critical point $\qbf^0=(s_1^+,s_1^-,\dots, s_n^+,s_n^-)$.
Indeed, consider the function
\begin{equation}
\label{eq:Psi}
\Psi(\qbf,\xbf)=\sum_{j=1}^n  (F_j(q_j^-,q_{j+1}^+,0)+S(q_{j+1}^+,x_{j+1})+S(x_j,q_{j}^-)),
\end{equation}
where $q_j^\pm\in B_j^\pm$, $x_j\in D_j$ and $q_{n+1}^\pm=q_1^\pm$, $x_{n+1}=x_1$.

By Proposition \ref{prop:reflect}, for fixed $ \qbf=(q_1^+,q_1^-,\dots, q_n^+,q_n^-)\in\calB$, the function
$\xbf\to\Psi(\qbf,\xbf)$ has a nondegenerate critical point $\xbf=\xbf(\qbf)\in N^n$, $x_j=\xi(q_j^-,q_j^+)$ and by (\ref{eq:g0}),
the critical value is
\begin{eqnarray*}
\Psi(\qbf,\xbf(\qbf))&=&\Crit_{\xbf}\sum_{j=1}^n
\left(F_j(q_j^-,q_{j+1}^+,0)+R_0(q_{j}^+,x_{j},q_{j}^-)\right)\\
&=& \sum_{j=1}^n  (F_j(q_j^-,q_{j+1}^+,0)+R_0(q_j^+,q_{j}^-))=\Phi_0(\qbf).
\end{eqnarray*}
For $\qbf=\qbf^0$, we have $\xbf(\qbf^0)=\xbf^0$.
On the other hand, for fixed $\xbf$, the function $\qbf\to\Psi(\qbf,\xbf)$ has a nondegenerate critical
point $\qbf=\qbf(\xbf)$ of the form (\ref{eq:q}).
The critical value is
$$
\Psi(\qbf(\xbf),\xbf)=\calA^{(n)}_\kbf(\xbf)=\sum_{j=1}^n J(\gamma_j),
$$
where $\gamma_j$ is a trajectory of system $(M,H_0=E)$ joining $x_j,x_{j+1}\in N$ and crossing $\Sigma_\rho$
at the points $q_j^-,q_{j+1}^+$.
By the assumption,  the function $ \calA^{(n)}_\kbf(\xbf)$ has a nondegenerate critical point $\xbf^0$. Then
$(\qbf^0,\xbf^0)$, $\qbf^0=\qbf(\xbf^0)$, is a nondegenerate critical point of $\Psi$.
Hence  $\qbf^0$ is a nondegenerate critical point of $\Phi_0$.

By (\ref{eq:Jmu}),
$$
\Phi_\mu(\qbf)=\Phi_0(\qbf)+O(\mu\ln |\mu|).
$$
By the implicit function theorem, for small $\mu\ne 0$, $\Phi_\mu(\qbf)$ has a nondegenerate
critical point  near $\qbf^0$
which defines a periodic orbit of the system $(M\setminus N,H_\mu)$
shadowing the chain $\gamma^0$.
\qed

\medskip

If there is a symmetry group $\Phi_\theta:M\to M$, then everything will be invariant
under $\Phi_\theta$, and we obtain a proof of Theorem \ref{thm:symm}.

\subsection{Proof of Theorem \ref{thm:hyp}}

It is   similar to the proof of Theorem \ref{thm:per}. We only need to check   uniformity.
Let $\Lambda\subset K^\Z\times N^\Z$ be a compact hyperbolic $\calT$-invariant set of admissible trajectories
of the DLS.

There exist a finite collection $\{\Omega_k\}_{k\in I}$   of compact\footnote{Topology on the set of collision orbits  $\gamma$  is defined  by
 reparametrizing   $\gamma$ proportionally to the arc length (in the metric $\|\,\cdot\,\|$),
 and using the topology in  $C^0([0,1],M)$.} sets of collision orbits $\gamma:[t_-,t_+]\to M$
such that collision chains $(\gamma_j)_{j\in\Z}$ corresponding to trajectories  $(\kbf,\xbf)\in\Lambda$
are concatenations of collision orbits $\gamma_j\in \Omega_{k_j}$.

 Collision orbits $\gamma\in\Omega_k$
join pairs of nonconjugate points $x_-(\gamma)\in N$ and $x_+(\gamma)\in N$ which  form compact sets
$$
X_k^\pm=\{x_\pm(\gamma):\gamma\in \Omega_k\}\subset N.
$$
Take open   sets $D_k^\pm\Subset\calD_E\cap  N$ such that   $X_k^\pm\Subset D_k^\pm$ for all $k\in I$.
Set $D=\bigcup (D_k^+\cup D_{k'}^-)$. Take sufficiently small $\rho>0$ and let $\Sigma_\rho$ be the corresponding set (\ref{eq:Sigma}).

For any $\gamma\in\Omega_k$ let $s_-(\gamma)$ and $s_+(\gamma)$  be the first and last intersection points with $\Sigma_\rho$.
By the definition of a collision orbit (\ref{eq:coll}), the angles between initial and final velocities $v_\pm(\gamma)$
 and $N$, and the collision speeds $\|v_\pm(\gamma)\|$ are bounded away from 0. Hence there exists $c>0$, independent of $\rho$ and $\gamma\in\Omega_k$, 
 such that $d(x_\pm(\gamma),s_\pm(\gamma))\le c\rho$.
Then
$$
 Y_k=\{(s_-(\gamma),s_+(\gamma))\in \Sigma_\rho^2:\gamma\in \Omega_k\}
 $$
 is a compact set contained in $\calD_E$. We can assume that for any $(q_-,q_+)\in Y_k$ there is unique $\gamma(q_-,q_+)\in\Omega_k$
 such that $s_\pm(\gamma)=q_\pm$.

There is $\mu_0>0$
such that for any $\mu\in (-\mu_0,\mu_0)$,  any collision orbit $\gamma\in \Omega_k$,  any pair of points $q_-,q_+$
in the set
$$
\hat Y_k=\{(q_-,q_+)\in \Sigma_\rho^2: d((q_-,q_+),Y_k)\le \eps\}.
$$
are joined  by a trajectory $\beta_\mu=\beta_\mu(q_-,q_+,k)$
of system $(M\setminus N,H_\mu=E)$ which is close to $\gamma(q_-,q_+)$.
This follows from compactness of $\Omega_k$ and nonconjugacy of  $s_\pm(\gamma)$ along $\gamma\in\Omega_k$.
Then the Maupertuis action
$$
F_k(q_-,q_+,\mu)=J_{\mu}(\beta_\mu)
$$
is a smooth function on $\hat Y_k$.

Every collision chain corresponding to a trajectory in  $\Lambda$ is a concatenation of collision orbits in $\{\Omega_k\}$. Let $\Pi_{k,k'}\subset \Omega_k\times\Omega_{k'}$, $(k,k')\in I^2$, be the compact set of all pairs $\gamma\in \Omega_k$, $\gamma'\in\Omega_{k'}$ of neighbor collision orbits in such concatenations. Let $s_+(\gamma)$, $s_-(\gamma')$
be the corresponding points in $\Sigma_\rho$. Set
$$
V_{kk'}=\{(s_+(\gamma),s_-(\gamma')): (\gamma,\gamma')\in \Pi_{k,k'}\}.
$$
There is a constant $c>0$ such that
$d(q_+,q_-)\le 2c\rho$  for all $(q_+,q_-)\in V_{k,k'}$.
If we take $C>2c$, then $V_{kk'}\subset P_{\rho}$, where $ P_{\rho}$
is the set (\ref{eq:Prho}) corresponding to $D$ and the constant $C$.

Let
$$
s_\pm(\gamma)=f(x_\pm(\gamma),u_\pm(\gamma)),\qquad x_\pm(\gamma)\in D_k^\pm,\quad |u_\pm(\gamma)|=\rho.
$$
By the  jump condition and compactness of $\Lambda$, if $\delta>0$ is small enough,
$$
|u_+(\gamma)+u_-(\gamma')|\ge \delta\rho\quad\mbox{for all}\; (\gamma,\gamma')\in \Pi_{k,k'}.
$$
Let $Q_{\rho}\subset P_{\rho}$ be the set (\ref{eq:Qrho})  corresponding to  $D$ and the constants $C,\delta>0$.
Then  $V_{kk'}\subset Q_\rho$.  There exist $\eps>0$ such that  $d((q_+,q_-),V_{kk'})<\delta$ implies $(q_+,q_-)\in Q_\rho$.

\begin{rem}
If also the no straight reflection condition holds, then $|u_+(\gamma)-u_-(\gamma')|\ge \delta\rho$, and so
$(s_+(\gamma),s_-(\gamma'))\in \hat Q_{\rho}$. Then  $V_{kk'}\subset \hat Q_\rho$.
\end{rem}

We take $\rho>0$ so small that Theorem \ref{thm:connect}   holds in $Q_{\rho}$.
Then there is  $\mu_0>0$ such that for any $\mu\in  I_{\mu_0}$, the points $(q_+,q_-)\in Q_{\rho}$
can be joined by a trajectory
$\alpha_\mu=\alpha_\mu(q_+,q_-)$ with action $J_\mu(\alpha_\mu)=R_\mu(q_+,q_-)$.

Let $(\kbf,\xbf^0)\in\Lambda$ be  a trajectory of the DLS,
and let $\gamma^0=(\gamma_j^0)_{j\in\Z}$, $(\gamma_j^0,\gamma_{j+1}^0)\in \Pi_{k_jk_{j+1}}$, be the corresponding  collision chain, where $\gamma_j^0$ connects the points $x_j^0,x_{j+1}^0\in D$
and intersects $\Sigma_\rho$ at the points $s_j^-=s_-(\gamma_j)$ and $s_{j+1}^+=s_+(\gamma_j)$. Then $(s_j^-,s_{j+1}^+)\in \hat Y_{k_j}$ and $(s_j^+,s_{j}^-)\in Q_\rho$.

As in (\ref{eq:Phi}), consider the formal functional
$$
\Phi_\mu(\qbf)=\sum_{j\in\Z}  (F_{k_j}(q_j^-,q_{j+1}^+,\mu)+R_\mu(q_j^+,q_{j}^-)),\qquad \qbf=(q_j^-,q_j^+)_{j\in\Z},
$$
where
$$
(q_j^-,q_{j+1}^+)\in \hat Y_{k_j},\quad (q_j^+,q_j^-)\in Q_{\rho}.
$$
The functional depends  on  $\kbf\in I^\Z$, but we do not show it in the notation.
As  in the proof of Theorem \ref{thm:per},  $\Phi_\mu(\qbf)$ is the action of an infinite concatenation $\hat\gamma$ of trajectories $\alpha_\mu(q_j^+,q_{j}^-)$ and $\beta_\mu(q_{j}^-,q_{j+1}^+,k_{j})$. 
The derivative $D\Phi_\mu(\qbf)=\Gamma_\mu(\qbf)$ makes sense, so critical points are well defined.
As in the proof of Theorem \ref{thm:per},  critical points of $\Phi_\mu$
correspond to trajectories of system $(M\setminus N,H_\mu=E)$
shadowing the collision chain $\gamma^0$.

Let us show that for $\mu=0$ the functional $\Phi_0$ has a uniformly nondegenerate critical point $\qbf^0=\qbf^0(\xbf^0) $:
\begin{equation}
\label{eq:Gamma0}
\Gamma_0(\qbf^0)=0,\qquad \|D\Gamma_0(\qbf^0)^{-1}\|_\infty\le C_2,\qquad \Gamma_0(\qbf)=D\Phi_0(\qbf),
\end{equation}
with $C_2=C_2(\Lambda)$ independent of the trajectory $(\kbf,\xbf^0)\in\Lambda$. 
The $l_\infty$ norm can be defined by using the Riemannian metric on $M$
to identify $D\Gamma_0(\qbf)$ with a linear operator on  an $l_\infty$  Banach space
$$
E_1\subset \prod_{j\in\Z}(T_{q_j^+}\Sigma_\rho\times T_{q_j^+}\Sigma_\rho)
$$
with the $l_\infty$ norm.
A simpler option is to use local coordinates.

We can  introduce  coordinate charts $O_j^\pm$ on $\Sigma_\rho$ containing the points $s_j^\pm$   by using e.g.\ the exponential maps
$\exp_{s_j^\pm}:T_{s_j^\pm}\Sigma_\rho\to \Sigma_\rho$.
Then we identify $O_j^\pm$ with a ball $\{q\in\R^{m-1}:|q-s_j^\pm|<\eps\}$, $m=\dim M$.
Then    $\qbf=(q_j^-,q_j^+)_{j\in\Z}$ is represented by a point in a ball
$$
Z_1=\{(\qbf^-,\qbf^+): \|\qbf^\pm-\sbf^\pm\|_\infty<\eps\}
$$
in the Banach space $E_1=l_\infty(\R^{m-1})$.
Thus $\Gamma_\mu$ is now a  map $\Gamma_\mu:Z_1\subset E_1\to E_1$, and $D\Gamma_\mu(\qbf)$ is a bounded operator in $E_1$.

Similarly we introduce local coordinates in a ball $D_j=B_\eps(x_j^0)\cap N\subset D_{k_j}^+\cap D_{k_{j+1}}^-$ by using e.g.\ the exponential map $\exp_{x_j^0}:T_{x_j^0}N\to N$.
Then we  identify $x_j\in D_j$ with a point in the ball\footnote{Recall that components of $N$ may have different dimensions.}  $\{x\in \R^{n_j}: |x-x_j^0|<\eps\}$.
Then for a trajectory $(\kbf,\xbf)$ we can regard $\xbf=(x_j)_{j\in\Z}$ as a point in  a ball $Z_2=\{\xbf:\|\xbf-\xbf^0\|_\infty<\eps\}$ in the $l_\infty$ Banach space
$$
E_2 =\{\xbf\in \prod_{j\in\Z}\R^{n(j)}: \|\xbf\|_\infty=\sup |x_j|<\infty\}.
$$

To show that $D\Gamma_0(\qbf^0):E_1\to E_1$ is invertible, as in  (\ref{eq:Psi}), consider the functional
$$
\Psi(\qbf,\xbf)=\sum_{j\in\Z}  (F_{k_j}(q_j^-,q_{j+1}^+,0)+S(q_{j+1}^+,x_{j+1})+S(x_j,q_{j}^-)),
$$
where   $(\qbf,\xbf)\in Z=Z_1\times Z_2$. The functional is formal, but its derivatives
$$
D_{\qbf}\Psi(\qbf,\xbf)=G_1(\qbf,\xbf),\quad  D_{\xbf}\Psi(\qbf,\xbf)=G_2(\qbf,\xbf),
$$
are well defined.
Then $G=(G_1,G_2)$ is a $C^1$ map from an open set $Z=Z_1\times Z_2\subset E=E_1\times E_2$ to  $E$ and
$$
\|G(\qbf,\xbf)\|_\infty\le C_3,\quad \|DG(\qbf,\xbf)\|_\infty\le C_3,\qquad (\qbf,\xbf)\in Z,
$$
where the constant $C_3=C_3(\Lambda)$ is independent of the trajectory $(\kbf,\xbf^0)$.

By Proposition \ref{prop:reflect}, the equation $G_2(\qbf,\xbf)=0$ has  a nondegenerate solution $\xbf(\qbf)$
such that
$$
G_2(\qbf,\xbf(\qbf))=0,\qquad \|D_\xbf G_2(\qbf,\xbf(\qbf))^{-1}\|_\infty\le C_4=C_4(\Lambda).
$$
Indeed, $\xbf(\qbf)=(x_j(\qbf))$ where $x_j(\qbf)=\xi(q_{j-1}^+,q_j^-)$ is  a nondegenerate critical point
 of the function
$x\to R(q_j^+,x,q_j^-)$. The operator $D_\xbf G_2$
is block diagonal, so nondegeneracy implies that the inverse is $l_\infty$ bounded.
We have $G_1(\qbf,\xbf(\qbf))=\Gamma_0(\qbf)$.

Now we use the following lemma \cite{Bol:DCDS} which is a  version of the Lyapunov--Schmidt reduction.

\begin{lem}\label{lem:nondeg_2}  
Let $E=E_1\times E_2$ be Banach spaces and let $G=(G_1,G_2):Z\to
E$ be a $C^1$ map of an open set $Z=Z_1\times Z_2\subset E$. Suppose that
there is $C>0$ such that
$$
\|DG(\qbf,\xbf)\|\le C,\quad \|(D_\xbf G_2(\qbf,\xbf))^{-1}\|\le C \quad\mbox{for
all}\quad (\qbf,\xbf)\in Z.
$$
Let $G(\qbf^0,\xbf^0)=0$ and let $\xbf=\xbf(\qbf)$, $\qbf\in  Z_1$,
be a solution of $G_2(\qbf,\xbf)=0$ such that $\xbf(\qbf^0)=\xbf^0$. Set
$\Gamma(\qbf)=G_1(\qbf,\xbf(\qbf))$.  Then $D\Gamma(\qbf)$ and $DG(\qbf,\xbf(\qbf))$ are invertible
simultaneously and there exists a constant $c=c(C)>0$ such that
$$
\|(D\Gamma(\qbf))^{-1}\|\le \|(DG(\qbf,\xbf(\qbf)))^{-1}\|\le
c(1+\|(D\Gamma(\qbf))^{-1}\|).
$$
\end{lem}

Thus (\ref{eq:Gamma0}) holds. We conclude that $\qbf^0$ is a uniformly $l_\infty$-nondegenerate critical point of $\Phi_0$
independently of a trajectory $(\kbf,\xbf^0)\in \Lambda$. By (\ref{eq:Jmu}),
$$
\|D\Gamma_\mu(\qbf)-D\Gamma_0(\qbf)\|_\infty\le C_5|\mu\ln|\mu||,\qquad C_5=C_5(\Lambda).
$$
The proof  of Theorem \ref{thm:hyp}  is completed by using a uniform version of the implicit function theorem.

\section{Regularization}

\label{sec:reg}

In this section we prove   Theorem \ref{thm:connect2} for $d=\codim  N\le 3$.
For $d\le2$ we use  the Levi-Civita regularization, and for $d=3$ the KS-regularization.
For $d\ge 4$ a different   method is needed.

Let  $f_\mu:T^\perp N\to M$ the exponential map
 corresponding to the Riemannian metric
$\|\,\cdot\,\|_\mu$. As in (\ref{eq:coord}), we assume that $D\Subset N$ is contractible
and choose an orthonormal basis $e_1(x),\dots,e_d(x)$ in $T_x^\perp N$ smoothly depending on $x\in D$.
The map
\begin{equation}
\label{eq:fmu}
D\times B_\rho\to U_\rho,\qquad q=f_\mu(x,u),\quad u=u_1e_1(x)+\dots+u_de_d(x),
\end{equation}
defines semigeodesic coordinates $x\in D$, $u\in B_\rho=\{ v\in \R^d:|v|<\rho\}$, in $U=U_\rho$. The Riemannian metric
$\|\,\cdot\,\|_\mu$ in $U$ has the form
\begin{equation}
\label{eq:metr}
\|\dot q\|_\mu^2=\langle \calA(x,u,\mu)\dot x,\dot x\rangle+|\dot u|^2 +\langle \calB(x,u,\mu)\dot u,\dot u\rangle+ \langle \calC(x,u,\mu)\dot x,\dot u\rangle,
\end{equation}
where $\calA$ is positive definite and\footnote{Here $O_k(u)$ means a function whose Taylor expansion with respect to $u$ starts with $k$-th order terms.}
$$
\calB(x,u,\mu)=O_2(u),\quad \calC(x,u,\mu)=O_1(u).
$$
By the properties of the exponential map, $d_\mu(q,N)=|u|$.

Let $y \in T_x^*D$, $v\in \R^d$, be the momenta conjugate to $x,u$, so that
$$
\langle p, dq\rangle=\langle y , dx\rangle+\langle v, du\rangle.
$$
Then
\begin{equation}
\label{eq:p_mu}
\|p\|_\mu^2=\langle  A(x,u,\mu)y ,y \rangle + |v|^2+\langle  B(x,u,\mu)v,v\rangle+ \langle  C(x,u,\mu)y ,v\rangle,
\end{equation}
where
$$
 B(x,u,\mu)=O_2(u),\quad  C(x,u,\mu)=O_1(u),\quad A(x,0,\mu)=\calA^{-1}(x,0,\mu).
$$

The gyroscopic 1-form  is
\begin{equation}
\label{eq:covector}
\langle w_\mu(q),dq\rangle=\langle a(x,u,\mu),dx\rangle +\langle b(x,u,\mu),du\rangle.
\end{equation}
Without loss of generality we may assume that
\begin{equation}
\label{eq:b}
b(x,u,\mu)=O(u).
\end{equation}

Indeed,
$$
\langle w_\mu(q),dq\rangle=\langle \tilde a(x,u,\mu),dx\rangle + \langle \tilde b(x,u,\mu),du\rangle+d\varphi(x,u,\mu),
$$
where
$$
 \tilde b=b(x,u,\mu)-b(x,0,\mu),\quad \varphi=\langle b(x,0,\mu),u\rangle,\quad \tilde a=a(x,u,\mu)-D_x\varphi(x,u,\mu).
$$
The differential $d\varphi$ can be dropped: it  does not affect trajectories $q(t)$ (only the corresponding momenta $p(t)$) 
since it changes only the boundary terms in the action functional (\ref{eq:J}).
The new coefficient $\tilde b$ satisfies   (\ref{eq:b}).

In the symplectic variables $x,y,u,v$ the  Hamiltonian (\ref{eq:Hmu}) has the form
\begin{eqnarray}
H_\mu(q ,p)&=&\frac12\Big(\langle  A(x,u,\mu)(y-a(x,u,\mu)) ,y-a(x,u,\mu) \rangle + |v-b(x,u,\mu)|^2\nonumber\\
&&+\langle  B(x,u,\mu)(v-b(x,u,\mu)),v-b(x,u,\mu)\rangle\nonumber\\
&&+ \langle C(x,u,\mu)(y-a(x,u,\mu) ),v-b(x,u,\mu)\rangle\Big)\nonumber\\
&&+  W(x,u,\mu)-\frac{\mu\phi(x,u,\mu)}{|u|}. \label{eq:normal}
\end{eqnarray}
Next we regularize the singularity at $u=0$.

\subsection{Codimension 2}

Let $d=2$. Then we   identify $\R^2=\C$  and   use the Levi-Civita  change of variables
$$
u=u(\xi)=\xi^2/2, \quad |u|=|\xi|^2/2, \quad du=\xi\,d\xi.
$$
In the real variables,
$$
  u(\xi)=\frac12 \Gamma(\xi)\xi  ,\qquad   \Gamma(\xi)=\left(
\begin{array}{cc}
\xi_1 &-\xi_2\\
\xi_2 &\xi_1
\end{array}
\right).
$$
The matrix $\Gamma$ is orthogonal:
\begin{equation}
\label{eq:Gamma}
  \Gamma^*(\xi)\Gamma(\xi)=\Gamma(\xi)\Gamma^*(\xi)=|\xi|^2 I_2.
\end{equation}
The square map  evidently satisfies $u(\xi_+)=u(\xi_-)$ iff $\xi_+=\pm \xi_-$ and
\begin{equation}
\label{eq:perp}
u(\xi_+)=-u(\xi_-)\quad\Leftrightarrow\quad \langle\xi_+,\xi_-\rangle=0,\quad |\xi_+|=|\xi_-|.
\end{equation}

Let $\eta$ be the   momentum  conjugate to $\xi$ so that
$$
\langle v,du\rangle=\langle v,\Gamma(\xi)d\xi\rangle=\langle \Gamma^*(\xi)v,d\xi\rangle= \langle\eta,d\xi\rangle.
$$
Thus
$$
\eta=\Gamma^*(\xi)v,\quad v= \frac{\Gamma(\xi) \eta}{|\xi|^2}, \quad |\eta|=|v||\xi|.
$$

\begin{rem}
In the complex notation, the formulas are much simpler: e.g.\ $\eta=\bar\xi v$.
But we need to write the transformation in the form which will work also for $d=3$.
\end{rem}

The   gyroscopic 1-form   is now
$$
\langle w_\mu(q),dq\rangle=\langle a(x,u(\xi),\mu),dx\rangle+ \langle \hat b(x,\xi,\mu), d\xi\rangle,
$$
where
$$
\hat b(x,\xi,\mu)=  \Gamma^*(\xi)b(x,u(\xi),\mu)=O_3(\xi).
$$

We have
$$
v-b(x,u,\mu)=\frac{\Gamma(\xi) \eta-|\xi|^2b(x,u(\xi),\mu)}{|\xi|^2}=\frac{\Gamma(\xi)( \eta- \hat b(x,\xi,\mu))}{|\xi|^2}.
$$
By (\ref{eq:p_mu}),
\begin{eqnarray*}
\|p-w_\mu(q)\|_\mu^2= \langle   A(x,u(\xi),\mu)(y-a(x,u(\xi),\mu)) ,y-a(x,u(\xi),\mu) \rangle\\
+\frac{|\eta-\hat b(x,\xi,\mu)|^2}{ |\xi|^2}
+\frac{\langle \hat B(x,\xi,\mu)(\eta- \hat b(x,\xi,\mu)),\eta- \hat b(x,\xi,\mu)\rangle}{|\xi|^4}\\
+\frac{\langle \hat C(x,\xi,\mu)(y -a(x,u(\xi),\mu)),\eta- \hat b(x,\xi,\mu)\rangle}{|\xi|^2}
\end{eqnarray*}
where
\begin{eqnarray}
\hat B(x,\xi,\mu) =\Gamma^*(\xi)B(x,u(\xi),\mu)\Gamma(\xi)=O_6(\xi),\label{eq:O6}\\
\hat C(x,\xi,\mu)=\Gamma^*(\xi)C(x,u(\xi),\mu)=O_3(\xi) .\label{eq:O3}
\end{eqnarray}

The equation $H_\mu=E$ takes the form
\begin{eqnarray*}
\left(\frac12 \langle  A(x,u(\xi),\mu)(y-a(x,u(\xi),\mu)) ,y-a(x,u(\xi),\mu) \rangle+W(x,u(\xi),\mu)-E\right)\frac{|\xi|^2}{2} \\
 +\frac{|\eta -\hat b(x,\xi,\mu)|^2}{2} +
\frac{\langle  \hat B(x,\xi,\mu)(\eta-\hat b(x,\xi,\mu), \eta- \hat b(x,\xi,\mu)\rangle}{2|\xi|^2} \\
+\frac12
\langle \hat C(x,\xi,\mu)(y-a(x,u(\xi),\mu) ,\eta-\hat b(x,\xi,\mu)\rangle =\mu\phi(x,u(\xi),\mu).
\label{eq:Hmu=E}
\end{eqnarray*}

Solving   for $\mu$ we obtain the regularized Hamiltonian
\begin{equation}
\label{eq:tildeH}
\mu=\calH(z,\xi,\eta)=\calH_2(z,\xi,\eta)+O_3(\xi,\eta) ,\qquad z=(x,y),
\end{equation}
where
$$
\calH_2(z,\xi,\eta)= \frac{(F_0(z)-E)|\xi|^2+|\eta|^2}{2\phi_0(x)},\qquad \phi_0(x)=\phi(x,0,0).
$$
Here
\begin{eqnarray*}
F_0(z)=H_0(x,y ,0,0)=\frac12\langle A_0(x)(y-a_0(x)),y-a_0(x)\rangle + W_0(x)
\end{eqnarray*}
is the Hamiltonian (\ref{eq:calM}) on $T^*N$ corresponding to the Lagrangian $L_0|_{TN}$.
Indeed,
$$
A_0(x)=A(x,0,0),   \quad a_0(x)=a(x,0,0),\quad W_0(x)=W(x,0,0).
$$ 

By (\ref{eq:O6}), the regularized Hamiltonian $\calH$ is at least of class $C^{3+\Lip}$, and the only source of low regularity is the term
$|\xi|^{-2}\hat B(x,\xi,\mu)=O_4(\xi)$. In applications to celestial mechanics,
$\hat B$ is divisible by $|\xi|^2$, so $\calH$ is real analytic.

Since  in the new symplectic  variables $x,y,\xi,\eta$ the level set $\{H_\mu=E\}$ becomes $\{\calH=\mu\}$, the symplectic map
$$
\psi(x,y,\xi,\eta)=(x,y,u(\xi),v(\xi,\eta))
$$
takes  solutions of the regularized Hamiltonian system on the level set $\{\calH=\mu\}$  to  solutions of the original Hamiltonian
system on the level set $\{H_\mu=E\}$ (with different time parametrization).

For fixed $z$, $\calH_2$ is a quadratic Hamiltonian with   eigenvalues
\begin{equation}
\label{eq:lambda}
\pm \lambda(z),\qquad \lambda(z)=\frac{\sqrt{E-F_0(z)}}{\phi_0(x)},
\end{equation}
each of multiplicity 2.
We see that $\xi=\eta=0$ is a critical manifold for $\calH$  and
 $\calM=\calM_E=\{(z,0,0):F_0(z)<E\}$ is a normally hyperbolic invariant manifold.

\medskip

Let $\tilde U=D\times B_r$, $r=\sqrt{2\rho}$ and let $\pi(x,\xi)=(x,u(\xi))$.
We proved the following semi global version of the Levi-Civita regularization:

\begin{thm}
\label{thm:Levi-Civita}
Let $D\Subset N\cap\calD_E$ be a domain  such that $T^\perp N|_D$ is trivial.
There exist a tubular neighborhood $U$ of $D$, a smooth map $\pi:\tilde U\to U$ and a $C^{3+\Lip}$ Hamiltonian $\calH$ on $T^*\tilde U$ such that:
\begin{itemize}
\item $\pi:\tilde U\setminus \tilde D\to U\setminus D$, $\tilde D=\pi^{-1}(D)$, is a double covering branched over $D$ and
$\pi:\tilde D\to  D$ a diffeomorphism;
\item
$\calH$ is invariant under the  sheet interchanging involution $\sigma:  \tilde U\to \tilde U$;
\item
$\pi$ takes trajectories of system $(\tilde U\setminus \tilde D,\calH=\mu)$    to trajectories
of system $(U\setminus D,H_\mu=E)$   (with changed time parametrization);
\item
The Hamiltonian system $(\tilde U,\calH)$ has  a $2(m-2)$-dimensional normally hyperbolic symplectic critical manifold
$\calM$ on the level $\calH=0$
with $2(m-2)$ zero eigenvalues and two semisimple real nonzero eigenvalues (\ref{eq:lambda}),
each of multiplicity 2.
\item
Trajectories of system $(\tilde U,\calH=0)$ asymptotic to $\calM_E$ are projected by $\pi$
to trajectories of the degenerate billiard   $(M, N,H=E)$ colliding with $N$.
\end{itemize}
\end{thm}

Since $\pi$ is a double covering, to each  orbit $\gamma:[0,\tau]\to U$ of the degenerate billiard colliding with $N$
at $x=\gamma(\tau)$,
there correspond 2 asymptotic trajectories $\gamma_{1,2}:[0,+\infty)\to \tilde U$, $\gamma_2=\sigma\gamma_1$,
of the regularized system with $\gamma_{1,2}(+\infty)= \pi^{-1}(x)$ and $\{\gamma_{1}(0),\gamma_2(0)\}=\pi^{-1}(\gamma(0))$. Similarly for an  orbit $\gamma:[\tau,0]\to U$ with $\gamma(\tau)\in N$.

\medskip

To prove  Theorem \ref{thm:connect2} we use a generalization of the Shilnikov lemma \cite{Shilnikov} for normally hyperbolic invariant manifolds of a Hamiltonian system. 
Let $\calM$ be a symplectic manifold with symplectic coordinates $z=(x,y)$. Consider a Hamiltonian system with Hamiltonian
\begin{equation}
\label{eq:calH}
\calH(z,\zeta)=\calH_2(z,\zeta)+O_3(\zeta),\qquad \calH_2=\frac12(a(z)|\eta|^2-b(z)|\xi|^2),
\end{equation}
where $z\in\calM$, $\zeta=(\xi,\eta)\in\R^{2d}$ and $a(z),b(z)>0$ for $z\in\calM$. Thus $(z_0,0,0)$ is a hyperbolic equilibrium with nonzero eigenvalues $\pm\lambda(z_0)$,
$\lambda(z_0)=\sqrt{a(z_0)b(z_0)}$. Its stable and unstable manifolds are given by
\begin{eqnarray*}
W^\pm(z_0)=\{(z,\xi,\eta): z=g_\pm(z_0,\xi),\; \eta=h_\pm(z_0,\xi)\},\qquad g_\pm(z_0,0)=z_0,\quad h_\pm(z_0,0)=0.
\end{eqnarray*}

Let $r>0$. Fix a compact set $\calK\subset\calM$ and   $\eps>0$ and denote
\begin{eqnarray}
\calQ_-=\{(z_0,\xi_-,\xi_+)\in \calK\times S_r^2: \langle \xi_-,\xi_+\rangle\ge \eps^2 r^2\}\\
\calQ_+=\{(z_0,\xi_-,\xi_+)\in \calK\times S_r^2: \langle \xi_-,\xi_+\rangle\le -\eps^2 r^2\}.
\end{eqnarray}

The next result is a  corollary  of Theorem 6 in \cite{Bol-Neg:RCD}.  

\begin{thm}
\label{thm:Shil} There exists $r>0$ and $\mu_0>0$ such that for
any
$$
(z_0,\xi_-,\xi_+,\mu)\in\calX= (\calQ_+\times (0,\mu_0))\cup  (\calQ_-\times (-\mu_0,0))
$$
\begin{itemize}
\item
There exists
$$
T\sim -\frac{1}{2\lambda(z_0)}\ln\left(\frac{-\mu}{\langle \xi_+,\xi_-\rangle}\right)
$$
and a solution
$$
\zeta(t)=(z(t),
 \xi(t),\eta(t))\in \calM\times  B_r\times \R^d ,\qquad t \in [- T, T],
$$
 with $\calH=\mu$  such that
\begin{equation}
\label{eq:BC0} z(0)=z_0,\quad \xi(T)=\xi_-,\quad
\xi(- T)= \xi_+.
\end{equation}
\item We have
\begin{eqnarray}
\xi(0)&=& \sqrt{\frac{-\mu}{2b(z_0)\langle \xi_+,\xi_-\rangle}}(\xi_++\xi_-)+O(\sqrt{|\mu|}r)+O(\mu),\label{eq:xi0}\\
\eta(0)&=& \sqrt{\frac{-\mu}{2a(z_0)\langle \xi_+,\xi_-\rangle}}(\xi_--\xi_+)+O(\sqrt{|\mu|}r)+O(\mu),\label{eq:eta0}\\
z_\pm&=&g_\pm(z_0,\xi_\pm)+O(\mu\ln|\mu|),\nonumber\\
\eta_\pm&=&h_\pm(z_0,\xi_\pm)+O(\mu)\nonumber.
\end{eqnarray}
 
\item
$\zeta$ smoothly depends on $(z_0,\xi_-,\xi_+,\mu)\in \calX$ and converges (as  a nonparametrized curve), as $\mu\to 0$,
to the concatenation of asymptotic trajectories $\zeta_+:[0,+\infty)\to \calM\times B_r\times\R^d$ and
$\zeta_-:(-\infty,0]\to \calM\times B_r\times\R^d$
in the stable and unstable manifolds $W^\pm(z_0)$ of $z_0\in\calM$:
\begin{eqnarray*}
\zeta_+(0)=(z_+,\xi_+,\eta_+)\in W^+(z_0),\quad\zeta_+(+\infty)=(z_0,0,0),\\
  \zeta_-(0)=(z_-,\xi_-,\eta_-)\in W^-(z_0),\quad \zeta_-(-\infty)=(z_0,0,0).
\end{eqnarray*}
\item
The Maupertuis action of $\zeta$ is a smooth function on $\calX$ and has the form
\begin{equation}
J(\zeta)=\int_\zeta y\,dx+\eta\,d\xi=J_-(z_0,\xi_-)+J_+(z_0,\xi_+)+O(\mu\ln |\mu|),
 \end{equation}
 where $J_\pm(z_0,\xi_\pm)$ are the actions of the asymptotic trajectories $\zeta_\pm$.
 \item
 If $\mu<0$, or $\mu>0$ and  $|\xi_+-\xi_-|\ge \eps r$, then $|\xi(t)|\ge c\sqrt{|\mu|}$ for $t \in [- T, T]$.
\end{itemize}
\end{thm}

\begin{rem}
In \cite{Bol-Neg:RCD} the proof was given for a smooth Hamiltonian $\calH$. This is enough for
applications in celestial mechanics. However, one can check that the proof works if $\calH\in C^{3+\Lip}$.
\end{rem}

Let us prove Theorem \ref{thm:connect2}  for $d=2$. For definiteness let $\mu>0$. Let $r=\sqrt{2\rho}$.
For given $u_\pm\in S_\rho$ with with $u_+\ne -u_-$ we can find $\xi_\pm\in S_r$
 such that $u_\pm=u(\xi_\pm)$ and $\langle\xi_+,\xi_-\rangle \ne 0$ by (\ref{eq:perp}).
Replacing  $\xi_+$ with $-\xi_+$ if necessary (they correspond to the same $u_+$) we may assume that $\langle\xi_+,\xi_-\rangle >0$.
We conclude  that there is $\eps>0$ such that if $|u_++u_-|\ge \delta\rho$,
we can find $\xi_\pm\in S_r$ with $u_\pm=u(\xi_\pm)$
such that $\langle \xi_-,\xi_+\rangle \ge \eps r^2$.  Then $\pi$ takes the trajectory
in Theorem \ref{thm:Shil} to a trajectory of system $(M\setminus N,H_\mu=E)$
satisfying the condition of  Proposition \ref{prop:reflect}.
\qed

\subsection{Codimension 3}

Let $d=3$. Then Theorem \ref{thm:Levi-Civita} is modified as follows:

\begin{thm}\label{thm:d=3}
There exist an $(m+1)$-dimensional   manifold $\tilde U$, a smooth group action
$\Phi_\theta:\tilde U\to\tilde U$, $\theta\in\T$,  a
smooth surjective map $\pi:\tilde U\to U$ commuting with $\Phi_\theta$,
and a  $\Phi_\theta$-invariant Hamiltonian $\calH\in C^{3+\Lip}$ on $T^*\tilde U$ such that:
\begin{itemize}
\item
The group action $\Phi_\theta$ is trivial on $\tilde D=\pi^{-1}(D)$
and free on $\tilde U\setminus \tilde D$. Thus $\pi:\tilde U\setminus \tilde D\to U\setminus D$
is a fiber bundle with fiber $\T$ and $\pi:\tilde D\to D$ is a diffeomorphism.
\item
Let  $G$ be the momentum integral
$$
G(q,p)=\langle X(q),p\rangle,\qquad  X(q)=D_\theta\big|_{\theta=0}\Phi_\theta(q)
$$
of system $(\tilde U,\calH)$ corresponding to the symmetry group $\Phi_\theta$. Then
$\pi$ takes trajectories of system $(\tilde U\setminus\tilde D,\calH=\mu)$ with  $G=0$ to trajectories
of system $(U,H_\mu=E)$.
\item System $(\tilde U,\calH)$ has  a $2(m-2)$-dimensional normally hyperbolic symplectic critical manifold
$\calM$ on the level  $\{\calH=0,G=0\}$.
 Every critical point $z\in \calM$ has $2(m-2)$ zero eigenvalues and two semisimple nonzero eigenvalues
(\ref{eq:lambda}), each of multiplicity 4.
\item
Trajectories asymptotic to $\calM$ are projected by $\pi$
to  trajectories colliding with $N$.
\end{itemize}
\end{thm}

Note that due to symmetry $\Phi_\theta$ to each  trajectory $\gamma:[0,\tau]\to U$ of the billiard colliding with $N$ 
there correspond a continuum of asymptotic orbits $\tilde\gamma:[0,+\infty)\to\tilde U$  of the regularized system with 
$\tilde\gamma(+\infty)=\pi^{-1}(\gamma(\tau))$ and $\tilde \gamma(0)\in\pi^{-1}(\gamma(0))$.

\medskip

\noindent{\it Proof.}
It is similar to the case $d=2$, only instead of the Levi-Civita
regularization we use the KS regularization \cite{K-S}.  The Hamiltonian still has the form  (\ref{eq:normal}),
but now $u,v\in\R^3$.
The square map $u:\R^2\to \R^2$  is replaced by the quadratic Hopf map\footnote{Quaternions provide a simpler formula for $u$.}
  $u:\R^4\to\R^3$
given by the Hurwitz matrix $\Gamma(\xi)$:
$$
u(\xi)=\frac 12\Gamma(\xi)\xi,\qquad \Gamma(\xi)=\left(
 \begin{array}{cccc}
 \xi_1&-\xi_2&-\xi_3&\xi_4\\
 \xi_2&\xi_1&-\xi_4 &-\xi_3\\
 \xi_3&\xi_4&\xi_1 &\xi_2
 \end{array}
 \right).
$$
 It has the following properties:
 \begin{itemize}
\item
 $|u(\xi)|=|\xi|^2/2$, $du(\xi)=\Gamma(\xi)\,d\xi$.
 \item
 We have
 $$
 \Gamma(e^{\theta J}\xi)=\Gamma(\xi)e^{\theta J} ,\qquad J=\left(
 \begin{array}{cc}
 0& -I_2\\
 I_2&0
 \end{array}
 \right).
 $$
 Thus $u(e^{\theta J} \xi)=u(\xi)$
 is invariant under the group $e^{\theta J}:\R^4\to\R^4$ generated by the vector field $J\xi$.
 \item
 $u(\xi_+)=u(\xi_-)$ iff $\xi_+=e^{\theta J}\xi_-$ for some $\theta$.
 \item
 $u(\xi_+)=-u(\xi_-)$ iff $|\xi_+|=|\xi_-|$ and $\langle \xi_+,\xi_-\rangle =0$, $\langle J\xi_+,\xi_-\rangle =0$.
 Equivalently, $\langle e^{\theta J}\xi_+,\xi_-\rangle\equiv 0$ for all $\theta$.
 \item
 $\Gamma(\xi)\Gamma^*(\xi)=|\xi|^2 I_3$.
 \item
 If $\langle J\xi,\eta\rangle=0$, then $\eta=\Gamma^*(\xi)v$ for a unique $v=v(\xi,\eta)\in\R^3$ given by
 $$
 v =\frac{\Gamma(\xi)\eta}{|\xi|^2},\qquad |\eta|=|\xi||v|,\qquad \langle v,du\rangle =\langle\eta,d\xi\rangle.
 $$
\end{itemize}

 Let
 $$
 Z=\{(\xi,\eta):\xi\ne 0,\;\langle J\xi,\eta\rangle=0\}
 $$
 and let $\tilde Z$ be the quotient of $Z$ under the group action
 $(\xi,\eta)\to (e^{\theta J}\xi,e^{\theta J}\eta)$. This is a symplectic manifold with a symplectic from derived from $d\eta\wedge d\xi$.
 The map
 $$
 (\xi,\eta)\in\tilde Z\to (u(\xi),v(\xi,\eta))\in (\R^3\setminus\{0\})\times \R^3
 $$
 is  invertible and it is  a symplectic diffeomorphism:
 $$
 \langle \eta,d\xi\rangle=\langle v(\xi,\eta),du(\xi)\rangle.
$$

We make a symplectic change of variables
$$
\psi: T^*D\times \tilde Z\to T^*D\times (\R^3\setminus \{0\})\times\R^3,\qquad \psi(x,y,\xi,\eta)=(x,y,u,v).
$$
 Define the regularized Hamiltonian $\calH(x,y,\xi,\eta)$ on $T^*D\times B_r\times \R^4$ by the same formula (\ref{eq:tildeH}), where $\hat B$ and $\hat C$
 are given by (\ref{eq:O6})--(\ref{eq:O3}).  It is easy to see that the Hamiltonian is invariant under the symplectic
 transformation $\Phi_\theta(z,\xi,\eta)=(z,e^{\theta J}\xi,e^{\theta J}\eta)$:
 $$
 \calH(z,e^{\theta J}\xi,e^{\theta J}\eta)=\calH(z,\xi,\eta),\qquad z=(\xi,\eta).
 $$
 Indeed,
 $$
 \hat b(x,e^{\theta J}\xi,\mu)=e^{\theta J}\hat b(x,\xi,\mu),\quad \hat B(x,e^{\theta J}\xi,\mu)=e^{\theta J}\hat B(x,\xi,\mu)e^{-\theta J}.
 $$
 Hence $\calH$ has the Noether integral $G=\langle J\xi,\eta\rangle$. On the zero level set $\{G=0\}$
 we have
 $$
 H_\mu(z,u(\xi),v(\xi,\eta))=E\quad \Leftrightarrow\quad \calH(z,\xi,\eta)=\mu.
 $$
By a standard result of the Hamiltonian reduction theory (see e.g.\ \cite{Arnold}), the map
 $\psi(z,\xi,\eta)=(z,u(\xi),v(\xi,\eta))$ takes trajectories of the regularized system  in $\{\calH=\mu\}\cap \{G=0\}$  to trajectories in $\{H_\mu=E\}$.
Theorem \ref{thm:d=3} is proved with $\tilde U=D\times B_r$, $r=\sqrt{2\rho}\}$ and $\pi(x,\xi)=(x,u(\xi))$.
\qed

\medskip

Let us prove Theorem \ref{thm:connect2}  for $d=3$. For definiteness let $\mu>0$. Let $r=\sqrt{2\rho}$.
For sufficiently small $\eps>0$ and given $u_\pm\in S_\rho$ with $|u_++u_-|\ge \delta\rho$, we need to find   $\xi_\pm\in S_r$ with $u_\pm=u(\xi_\pm)$
such that $\langle \xi_-,\xi_+\rangle \le -\eps r^2$. Then we can join $\xi_+$ and $\xi_-$  by a trajectory $\zeta(t)=(z(t),\xi(t),\eta(t))$
in  Theorem \ref{thm:Shil}. We will show that it is possible to choose $\xi_\pm$ in such a way that
this trajectory satisfies $G=\langle J\xi,\eta\rangle \equiv 0$.
Then $\pi$ takes the trajectory
in Theorem \ref{thm:Shil} to a trajectory of system $(U\setminus D,H_\mu=E)$
satisfying the conditions of  Theorem \ref{thm:connect2}.

Let us compute the value of $G$ along the trajectory $\zeta(t)$ in  Theorem \ref{thm:Shil}. By (\ref{eq:xi0})--(\ref{eq:eta0}),
$$
G=\langle J\xi(0),\eta(0)\rangle=-\frac{\mu\langle J\xi_+,\xi_-\rangle}{2\lambda(z_0)\langle \xi_+,\xi_-\rangle}+O(\mu r)=G(z_0,\xi_+,\xi_-,\mu).
$$

In the next computation we follow \cite{Bol-Mac:centers}. Suppose that $u_++u_-\ne 0$.
Then $u_\pm=u(\xi_\pm)$, where $s(\theta)=\langle e^{\theta J}\xi_+,\xi_-\rangle \not\equiv 0$.
Let $\theta_0$ be a maximum point of $s(\theta)$. Then
$$
s(\theta_0)=\langle e^{\theta_0 J}\xi_+,\xi_-\rangle>0,\quad s'(\theta_0)=\langle Je^{\theta_0 J}\xi_+,\xi_-\rangle=0,
$$
and the critical point is nondegenerate. By the implicit function theorem for small enough $r$ and $\mu_0$
there is $\theta$ near $\theta_0$ such that
$$
G(z_0,e^{\theta J}\xi_+,\xi_-,\mu)=0,\quad \langle e^{\theta J}\xi_+,\xi_-\rangle>0
$$
Then $\tilde \xi_+=e^{\theta J}\xi_+$ satisfies
$$
(\tilde \xi_+,\xi_-)\in\calQ_+,\quad G(z_0,\xi_+,\xi_-,\mu)=0.
$$
Now the trajectory in Theorem \ref{thm:Shil} corresponding to $z_0,\tilde\xi_+,\xi_-$  
 is projected by $\pi$ to a trajectory of system $(\tilde U,H_\mu=E)$ satisfying the conditions of  Theorem \ref{thm:connect2}.
\qed

 \end{document}